\input amstex.tex
\input amsppt.sty

\input boxedeps.tex

\HideDisplacementBoxes

\magnification 1000
\pageheight{9truein} \pagewidth{6.5 truein}

\def\bi{\it}

\nologo\TagsOnRight
\NoBlackBoxes

\newif\ifcontains    
\def\contains#1#2{\def\test##1#1##2##3\endtest{
\ifx##2\ZZZ%
\containsfalse\else\containstrue\fi}\test#2#1\ZZZ
\endtest}

\def\cite#1{%
 \def\commasplit##1,##2\endcommasplit{[$
\bold{##1}$,##2]}%
 \contains,{#1}%
 \ifcontains \commasplit#1\endcommasplit%
 \else [{\bf#1}]\fi}

\def\qed{\ifhmode\unskip\nobreak\fi\quad\hfill
\raise.32em\hbox{\boxed{}}} 

\topmatter\title
A Schl{\"a}fli-type formula\\
for convex cores of hyperbolic 3--manifolds
\endtitle
\author
Francis Bonahon
\endauthor
\affil
University of Southern California
\endaffil
\address
Department of Mathematics, 
University of Southern California, 
Los Angeles CA 90089-1113, U.S.A.
\endaddress
\dedicatory
\hfill Dedicated to D.B.A.~Epstein, on his 60th 
birthday.
\enddedicatory
\email fbonahon\@math.usc.edu \endemail

\thanks
This research was partially supported by N.S.F. 
grants DMS-9001895, DMS-9201466 and DMS-9504282.
\endthanks
\endtopmatter
\document

\rightheadtext{A Schl{\"a}fli-type formula}

\nobreak\par\nobreak 	Let $M$ be a (connected) 
hyperbolic 3--manifold, namely a complete 
Riemannian manifold of dimension 3 and of 
constant sectional curvature $-1$, with finitely 
generated fundamental group. A fundamental subset 
of $M$ is its {\bi convex core \/} $C_{M}$, which 
is the smallest non-empty convex subset of $M$. 
The condition that the volume of $C_{M}$ is 
finite is open in the space of hyperbolic metrics 
on $M$, provided we restrict attention to 
cusp-respecting deformations. In this paper, we 
give a formula which, for a cusp-preserving 
variation of the hyperbolic metric of $M$, 
expresses the variation of the volume of the 
convex core $C_{M}$ in terms of the variation of 
the bending measure of its boundary.
\nobreak\par\nobreak 	This formula is analogous 
to the Schl{\"a}fli formula for the volume of an 
$n$--dimensional hyperbolic polyhedron $P$; see 
\cite{Sc1}\cite{Kne}\cite{AVS} and {\S }1. If the 
metric of $P$ varies, the Schl{\"a}fli formula 
expresses the variation of the volume of $P$ in 
terms of the variation of the dihedral angles of 
$P$ along the $\left ( n-2\right )$--faces of its 
boundary and of the $\left ( n-2\right 
)$--volumes of these faces. 
\nobreak\par\nobreak 	The analogy stems from the 
fact that the boundary $\partial C_{M}$ of 
$C_{M}$ is almost polyhedral, in the sense that 
it is totally geodesic almost everywhere. 
However, the {\bi pleating locus \/}, where 
$\partial C_{M}$ is not totally geodesic, is not 
a finite collection of edges any more. Typically, 
it will consist of uncountably many infinite 
geodesics. In addition, the topology of this 
pleating locus can drastically change as we vary 
the metric of $M$. So the situation is much more 
complex.
\nobreak\par\nobreak 	The path metric induced on 
the surface $\partial C_{M}$ by the metric of $M$ 
is hyperbolic with finite area. On this 
hyperbolic surface, the pleating locus $\lambda $ 
forms a {\bi compact geodesic lamination \/}, 
namely is compact and is the union of disjoint 
simple geodesics. The surface $\partial C_{M}$ is 
bent along $\lambda $, and the amount of this 
bending can be measured, not by dihedral angles 
any more, but by a transverse measure for 
$\lambda $. Endowing $\lambda $ with this 
transverse measure, we get a measured lamination 
$b$, called the {\bi bending measured lamination 
\/} of $M$; see \cite{Thu}\cite{EpM}. 
\nobreak\par\nobreak 	Let $M$ be a hyperbolic 
3--manifold which is {\bi geometrically finite 
\/}, namely such that the convex core $C_{M}$ has 
finite volume and such that the fundamental group 
$\pi _{1}\left ( M\right )$ is finitely 
generated. Consider a deformation of $M$, namely 
a differentiable 1--parameter family of 
hyperbolic manifolds $M_{t}$, $t\in \left [0, 
\varepsilon \right [$, such that $M_{0}=M$; when 
$M$ has cusps, we also require that the cusps of 
each $M_{t}$ precisely correspond to the cusps of 
$M$. Then, $M_{t}$ is also geometrically finite 
for $t$ small enough \cite{Mar}. We showed in 
\cite{Bo4} that, if  $b_{t}$ is the bending 
measured lamination of $M_{t}$, then the family 
$b_{t}$, $t\in \left [0, \varepsilon \right [$, 
admits a tangent vector $\dot b_{0}$ at $t=0$, in 
the piecewise linear manifold ${\Cal M}{\Cal 
L}\left ( \partial C_{M}\right )$ of all measured 
geodesic laminations on $\partial C_{M}$. In 
addition, in \cite{Bo1}\cite{Bo2}, we showed that 
such a tangent vector can be geometrically 
interpreted as a geodesic lamination endowed with 
a certain type of transverse distribution, called 
a transverse H{\"o}lder distribution.
\nobreak\par\nobreak 	On a hyperbolic surface, a 
geodesic lamination with transverse distribution 
$a$ admits a certain length \cite{Bo2}\cite{Bo3}. 
This length is designed so that it varies 
continuously with $a$ and coincides with the 
usual length when $a$ consists of a simple closed 
geodesic endowed with the Dirac transverse 
distribution. In particular, on the hyperbolic 
surface $\partial C_{M_{0}}$, we can consider the 
length $l_{0}\bigl ( \dot b_{0}\bigr ) $ of the 
tangent vector $\dot b_{0}$. 

\proclaim {Main Theorem }   With the above data, 
the volume $V_{t}$ of the convex core $C_{M_{t}}$ 
admits a right derivative $\dot V_{0}$ at $t=0$, 
and
$$\dot V_{0} = {\textstyle {1\over 2}} l_{0}\bigl 
( \dot b_{0}\bigr ) $$
where $l_{0}\bigl ( \dot b_{0}\bigr ) $ is the 
length of the vector $\dot b_{0}$ tangent to the 
family of bending measured laminations $b_{t}$.  
\endproclaim    

\nobreak\par\nobreak 	In the case of a 
differentiable deformation $M_{t}$, $t\in \left 
]-\varepsilon ,\varepsilon \right [$, the right 
and left derivatives of the volume of $C_{M_{t}}$ 
may not necessarily agree at $t=0$, as shown for 
instance by the example of \cite{Bo4, {\S }6}. 
\nobreak\par\nobreak 	An application of this 
theorem is the following corollary, proved in {\S 
}4.

\proclaim {Corollary }   Given a geometrically 
finite hyperbolic $3$--manifold $M$, consider the 
volumes of the convex cores of the 
cusp-preserving deformations of $M$. If the 
boundary of the convex core $C_{M}$ is totally 
geodesic, then $M$ corresponds to a local minimum 
of this volume function. \endproclaim    

\nobreak\par\nobreak 	We prove the Main Theorem 
in two steps. One step, proved in {\S }3, is to 
show that the volume of the convex core of 
$M_{t}$ has the same right derivative at $t=0$ as 
the volume enclosed in $M_{t}$ by a pleated 
surface whose pleating locus is constant and 
contains the pleating locus of $\partial 
C_{M_{0}}$. This step heavily relies on the 
arguments on \cite{Bo4}. The other step, proved 
in {\S }2, is devoted to a Schl{\"a}fli formula 
for the volume enclosed by a pleated surface 
whose pleating locus is constant. This simpler 
formula is proved by cutting the enclosed volume 
into small pieces and applying the usual 
Schl{\"a}fli formula to the pieces. The formal 
aspects of this part of the proof are relatively 
natural. However, much care is needed to justify 
these formal arguments, owing to the subtleties 
of the convergence of transverse distributions 
and to the fact that one has to estimate 
derivatives of dihedral angles, and not just 
dihedral angles. The vanishing of the 
contributions of internal edges is also 
non-trivial.
\medskip\nobreak 

\nobreak\par\nobreak 	This article was written 
while the author was holding visiting positions 
at the Centre \'Emile Borel, the Institut des 
Hautes \'Etudes Scientifiques and the California 
Institute of Technology. He would like to thank 
these institutions for their productive 
hospitality. He is also grateful to Steve Carlip 
for asking the question which originally 
motivated his interest in these problems, and to 
Alberto Candel for help in the proof of Lemma~9.

\head {\S }1. The Schl{\"a}fli formula for 
hyperbolic cycles in 3--manifolds\endhead    

\nobreak\par\nobreak 	The classical Schl{\"a}fli 
Formula is a crucial tool in the proof of Main 
Theorem. Although it holds in any dimension and 
in any space of non-zero constant curvature, we 
state it here only for hyperbolic 3--dimensional 
geometry since this is the only case which we 
will use. See \cite{Sc1}\cite{Kne}\cite{AVS, 
Chap. 7, {\S }2.2} for a proof.
\nobreak\par\nobreak 	Consider a differentiable 
1--parameter family of hyperbolic polyhedra 
$P_{t}$, $t\in \left [0, \varepsilon \right [$, 
in hyperbolic 3--space ${\Bbb H}^{3}$. This means 
that the polyhedra $P_{t}$ all have the same 
combinatorial type, that their faces and edges 
are totally geodesic in ${\Bbb H}^{3}$, and that 
their vertices vary differentiably with $t$.

\proclaim {Theorem 1 {\rm (Schl{\"a}fli Formula)} 
}   Let $P_{t}$, $t\in \left [0, \varepsilon 
\right [$, be a differentiable $1$--parameter 
family of polyhedra in ${\Bbb H}^{3}$. Then the 
right derivative of the volume $V_{t}$ of $P_{t}$ 
at $t=0$ is
$$\dot V_{0} = {\textstyle {1\over 2}} \sum 
_{e\roman {\thinspace edge\thinspace of\thinspace 
}P_{0}}  l_{0}\left ( e\right ) \dot b_{0}\left ( 
e\right )            \tag  1$$
where $l_{0}\left ( e\right )$ denotes the length 
of the edge $e$ in $P_{0}$, and where $\dot 
b_{0}\left ( e\right )$ is the right derivative 
at $t=0$ of the external dihedral angle 
$b_{t}\left ( e\right )$ of $P_{t}$ along this 
same edge.   \qed   \endproclaim

\nobreak\par\nobreak 	Here, the external dihedral 
angle $b_{t}\left ( e\right )$ is $\pi $ minus 
the internal dihedral angle of $P_{t}$ at $e$. In 
particular, the external dihedral angle is equal 
to 0 when the boundary of $P_{t}$ is flat at $e$, 
and is equal to $\pi $ when the two faces that 
are adjacent to $e$ locally coincide near $e$. 
\nobreak\par\nobreak 	There is a convenient 
notation, which already appeared in the 
Introduction and in the above statement, which we 
will consistently use throughout the paper: When 
a quantity $A_{t}$ depends on $t$, we will denote 
by $\dot A_{t_{0}}$ its right derivative with 
respect to $t$ at $t=t_{0}$. For instance, $\dot 
b_{0}\left ( e\right )$ is the right derivative 
of $b_{t}\left ( e\right )$ at $t=0$. 
\nobreak\par\nobreak 	Also, before going any 
further, we should observe that it suffices to 
prove the Main Theorem for orientable manifolds. 
Indeed, passing to the orientation cover 
multiplies each side of the equality by 2. 
Consequently, we will henceforth assume that all 
manifolds considered are orientable, and often 
oriented.
\nobreak\par\nobreak 	We can give a homological 
flavor to the formula of Theorem~1 in the 
following way. Let $M_{t}$, $t\in \left 
[0,\varepsilon \right [$, be a differentiable 
1--parameter family of connected oriented 
hyperbolic 3--manifolds. By definition, the 
differentiability condition means that there is a 
fixed group $\Gamma $ and orientation-preserving 
discrete faithful representations $\rho 
_{t}:\Gamma \rightarrow \roman {Isom}^{+}\left ( 
{\Bbb H}^{3}\right )$ into the isometry group of 
${\Bbb H}^{3}$ such that each $M_{t}$ is 
isometric to ${\Bbb H}^{3}/\rho _{t}\left ( 
\Gamma \right )$ and such that $\rho _{t}\left ( 
\gamma \right )$ depends differentiably on $t$ 
for every $\gamma \in \Gamma $. Fix a compact 
triangulated surface $S$ without boundary, not 
necessarily connected, and consider {\bi 
polyhedral maps \/} $f_{t}:S\rightarrow M_{t}$, 
namely continuous maps whose restriction to each 
edge and face of the triangulation of $S$ is a 
totally geodesic immersion. In addition, we 
require these maps  $f_{t}:S\rightarrow 
M_{t}={\Bbb H}^{3}/\rho _{t}\left ( \Gamma \right 
)$ to depend differentiably on $t$ in the sense 
that, if we lift them to maps $\widetilde  
f_{t}:\widetilde  S\rightarrow {\Bbb H}^{3}$ 
defined on the universal covering $\widetilde  S$ 
of $S$, the images of the vertices of $\widetilde 
 S$ under $\widetilde  f_{t}$ depend 
differentiably on $t$. 

\proclaim {Corollary 2 }   Given a triangulated 
compact oriented surface $S$, let 
$f_{t}:S\rightarrow M_{t}$, $t\in \left 
[0,\varepsilon \right [$, be a differentiable 
$1$--parameter family of polyhedral maps from $S$ 
to oriented hyperbolic $3$--manifolds $M_{t}$. 
Assume that the $f_{t}$ are homologous to 0 in 
$M_{t}$ and that the $M_{t}$ are non-compact, so 
that the volume $V_{t}$ of a $3$--chain bounding 
$f_{t}$ in $M_{t}$ is well defined.  Then the 
right derivative of the volume $V_{t}$ at $t=0$ 
is
$$\dot V_{0} = {\textstyle {1\over 2}} \sum 
_{e\roman {\thinspace edge\thinspace of\thinspace 
}S}  l_{0}\left ( e\right ) \dot b_{0}\left ( 
e\right )  $$
where, for each edge $e$ of $S$, $l_{0}\left ( 
e\right )$ denotes the length of $f_{0}\left ( 
e\right )$ and $b_{t}\left ( e\right )$ is the 
external angle between the two faces of 
$f_{t}\left ( S\right )$ meeting along 
$f_{t}\left ( e\right )$.  \endproclaim    
\demo {Proof}   Since $f_{0}$ is homologous to 0, 
we can extend it to a map $f_{0}:\Sigma 
\rightarrow M_{0}$ where $\Sigma $ is a 
simplicial complex with boundary $S$. We can 
choose this extension to be polyhedral. Then, 
since the hyperbolic manifolds $M_{t}$ and the 
maps $f_{t}:S\rightarrow M_{t}$ depend 
differentiably on $t$, we can easily extend them 
to a differentiable 1--parameter family of 
polyhedral maps $f_{t}:\Sigma \rightarrow M_{t}$ 
for $t$ small enough. 
\nobreak\par\nobreak 	Let $P_{1}$, \dots , 
$P_{n}$ be the 3--simplices of $\Sigma $. Apply 
Theorem~1 to each hyperbolic 3--simplex $f_{t 
\vert P_{i}}:P_{i}\rightarrow {\Bbb H}^{3}$. Note 
that the volume $V_{t}^{i}$ of this simplex is 
negative when $f_{t \vert P_{i}}$ is orientation 
reversing. Then,
$$\align
\dot V_{0} &=  \sum _{i=1}^{n} \dot V_{0}^{i} = 
{\textstyle {1\over 2}}  \sum _{i=1}^{n} \sum 
_{e\roman {\thinspace edge\thinspace of\thinspace 
}P_{i}} l_{0}\left ( e\right ) \dot 
b_{0}^{i}\left ( e\right )  \\
		&= {\textstyle {1\over 2}} \sum _{e\roman 
{\thinspace edge\thinspace of\thinspace }\Sigma } 
l_{0}\left ( e\right )\sum _{P_{i}\roman 
{\thinspace containing\thinspace }e}  \dot 
b_{0}^{i}\left ( e\right ) 
\endalign$$
where $b_{t}^{i}\left ( e\right )$ is the 
external angle between the two faces of $f_{t 
\vert P_{i}}$ meeting along the edge $f_{t \vert 
e}$ (counted negative if $f_{t \vert P_{i}}$ is 
orientation reversing).
\nobreak\par\nobreak 	For every edge $e$ of 
$\Sigma $ that is not in $S$, 
$$\sum _{P_{i}\roman {\thinspace 
containing\thinspace }e} b_{t}^{i}\left ( e\right 
) \equiv  0 \roman {\thinspace mod\thinspace } 
2\pi $$
and it follows that the corresponding derivative 
is equal to 0. On the other hand, for every edge 
$e$ of $S$, 
$$\sum _{P_{i}\roman {\thinspace 
containing\thinspace }e} b_{t}^{i}\left ( e\right 
) \equiv  b_{t}\left ( e\right ) \roman 
{\thinspace mod\thinspace } 2\pi .$$
The formula of Corollary~2 immediately follows. 
\qed\enddemo

\head {\S }2. Cycles bounded by pleated surfaces 
in hyperbolic 3--manifolds\endhead    
\nobreak\par\nobreak 	Let $M_{t}$, $t\in \left 
[0,\varepsilon \right [$, be a differentiable 
1--parameter family of hyperbolic 3--manifolds, 
associated to the representations $\rho 
_{t}:\Gamma \rightarrow \roman {Isom}^{+}\left ( 
{\Bbb H}^{3}\right )$. We require that this 
deformation of $M_{0}$ is {\bi cusp-preserving 
\/}, in the sense that every element of $\Gamma $ 
which is sent to a parabolic element by $\rho 
_{0}$ is also sent to a parabolic element by each 
$\rho _{t}$. We also assume that $M_{0}$ is {\bi 
geometrically finite \/}, namely that $\Gamma $ 
is finitely generated and that the convex core 
$C_{M_{0}}$ has finite volume. Then, the same 
also holds for every $C_{M_{t}}$ with $t$ small 
enough \cite{Mar, {\S }9}. In addition, the 
topology of $M_{t}$ and $\partial C_{M_{t}}$ 
remains constant for $t$ small enough, provided 
we use the following convention: When $C_{M_{0}}$ 
is 2--dimensional, namely when the group $\rho 
_{0}\left ( \Gamma \right )\subset \roman 
{Isom}^{+}\left ( {\Bbb H}^{3}\right )$ respects 
a hyperbolic plane in ${\Bbb H}^{3}$, we define 
$\partial C_{M_{0}}$ as the orientation covering 
of $C_{M_{0}}$, namely as the two sides of 
$C_{M_{0}}$ in $M_{0}$ (in contrast to the 
topological convention for which $\partial 
C_{M_{0}}$ should be equal to $C_{M_{0}}$ in this 
case). 
\nobreak\par\nobreak 	We want to compute the 
variation of the volume of the convex core 
$C_{M_{t}}$, namely of the part of $M_{t}$ 
bounded by the pleated surface $\partial 
C_{M_{t}}$. As $t$ varies, the pleating locus of 
$\partial C_{M_{t}}$ usually changes, which is a 
source of technical difficulties. In this 
section, as a first step towards our goal, we 
consider a simpler situation by substituting to 
$\partial C_{M_{t}}$ a pleated surface in $M_{t}$ 
whose pleating locus is independent of $t$, and 
by considering the variation of the volume 
bounded by this pleated surface. The fact that 
the pleating locus is constant makes the 
situation reminiscent of that of the classical 
Schl{\"a}fli formula of Theorem~1. 
\nobreak\par\nobreak 	Let $S$ be an oriented 
surface of finite topological type, without 
boundary but possibly infinite and not 
necessarily connected. Consider a pleated surface 
$f_{0}:S\rightarrow M_{0}$. We refer to 
\cite{Thu}\cite{CEG, {\S }5} for basic facts 
about pleated surfaces. In particular, an 
important convention is that $f_{0}$ is proper, 
and sends each end of $S$ to a cusp of $M_{0}$; 
as a consequence, the hyperbolic metric $m_{0}$ 
of $S$ obtained from the metric of $M_{0}$ by 
pull back under $f_{0}$ has finite volume, and 
each end of $S$ corresponds to a cusp of this 
metric. Although this is not absolutely necessary 
(see Remark~2 at the end of this section), we 
also require that $f_{0}$ is totally geodesic 
near the ends of $S$.
\nobreak\par\nobreak 	Let $\lambda $ be a {\bi 
pleating locus \/} for $f_{0}$, namely a compact 
geodesic lamination in $S$ such that $f_{0}$ 
sends each leaf of $\lambda $ to a geodesic of 
$M_{0}$ and such that $f_{0}$ is a totally 
geodesic immersion on $S-\lambda $. Such a 
pleating locus may not be unique; an extreme 
example occurs when $f_{0}$ is totally geodesic, 
in which case every compact geodesic lamination 
in $S$ is a pleating locus for $f_{0}$. 
Increasing $\lambda $ without loss of generality 
(compare \cite{CEG, {\S }4}), we can assume that 
$\lambda $ is a {\bi maximal \/} among compact 
geodesic laminations, namely that every component 
of $S-\lambda $ is either an ideal triangle, 
bounded by 3 leaves of $\lambda $, or an open 
annulus bounded on one side by a leaf of $\lambda 
$ with one spike and leading to a cusp on the 
other side. Then, for every $t$ small enough, 
there is a unique pleated surface 
$f_{t}:S\rightarrow M_{t}$ with pleating locus 
$\lambda $ such that, for every leaf $g$ of 
$\lambda $, $f_{t}\left ( g\right )$ is the 
geodesic of $M_{t}$ that is asymptotic to the 
image of $f_{0}\left ( g\right )$ under the 
quasi-isometric homeomorphism $\varphi 
_{t}:M_{0}\rightarrow M_{t}$; see 
\cite{Thu}\cite{CEG, {\S }5.3}. 
\nobreak\par\nobreak 	In \cite{Bo3}, we describe 
the local geometry of the pleated surface $f_{t}$ 
by the hyperbolic metric it induces on $S$, as 
well as a  {\bi bending transverse cocycle \/} 
$b_{t}\in {\Cal H}\left ( \lambda ;{\Bbb R}/2\pi 
{\Bbb Z}\right )$ for the geodesic lamination 
$\lambda $, valued in ${\Bbb R}/2\pi {\Bbb Z}$, 
which measures the bending of the pleated surface 
$f_{t}$. This bending transverse cocycle 
associates a number $b_{t}\left ( k\right )\in 
{\Bbb R}/2\pi {\Bbb Z}$ to each arc $k$ 
transverse to $\lambda $, which measures the 
bending of the pleated surface $f_{t}$ along the 
leaves of $\lambda $ meeting $k$. This 
$b_{t}\left ( k\right )$ is invariant under 
homotopy of $k$ respecting $\lambda $, and 
behaves additively if we split $k$ into two 
subarcs. We also prove in \cite{Bo3} that this 
bending cocycle depends differentiably of the 
representation $\rho _{t}:\Gamma \rightarrow 
\roman {Isom}^{+}\left ( {\Bbb H}^{3}\right )$ 
associated to $M_{t}$. In particular, there 
exists a right derivative $\dot b_{0}\in {\Cal 
H}\left ( \lambda ;{\Bbb R}\right )$, which is an 
${\Bbb R}$--valued transverse cocycle for 
$\lambda $. 
\nobreak\par\nobreak 	In \cite{Bo1}, we showed 
that every real-valued transverse cocycle $b\in 
{\Cal H}\left ( \lambda ;{\Bbb R}\right )$ 
defines a transverse distribution for $\lambda $. 
In particular, given a finite area hyperbolic 
metric $m$ on $S$, we can define a length 
$l_{m}\left ( b\right )\in {\Bbb R}$ by, first 
making $\lambda $ an $m$--geodesic lamination, 
and then locally integrating with respect to the 
transverse distribution associated to $b$ the 
1--dimensional Lebesgue measure along the leaves 
of $\lambda $; see \cite{Bo2}\cite{Bo3}.
\nobreak\par\nobreak 	Finally, we assume that the 
pleated surfaces $f_{t}$ separate $M_{t}$. When 
$S$ is non-compact, this means that the locally 
finite 2--chain defined by $f_{t}$ bounds a 
locally finite 3--chain. We also require that 
this 3--chain has finite volume $V_{t}$. If 
$M_{t}$ has infinite volume, this finite volume 
$V_{t}$ is uniquely defined, namely is 
independent to the finite volume 3--chain 
bounding $f_{t}$. If $M_{t}$ has finite volume, 
$V_{t}$ is defined only modulo the volume of 
$M_{t}$; however Mostow's Rigidity Theorem 
implies that the $M_{t}$ and $f_{t}$ are 
independent of $t$ up to isometry, so that the 
theorem below is trivial in this case. 

\proclaim {Theorem~3 }   Given an oriented 
surface $S$ of finite topological type, let 
$f_{t}:S\rightarrow M_{t}$, $t\in \left 
[0,\varepsilon \right [$, be a differentiable 
$1$--parameter family of pleated surfaces in 
oriented geometrically finite hyperbolic 
$3$--manifolds $M_{t}$, with pleating locus a 
fixed compact geodesic lamination $\lambda $ in 
$S$. Assume that $f_{t}$ bounds a finite volume 
(locally finite) $3$--chain in $M_{t}$ and that 
$M_{t}$ has infinite volume, so that we can 
consider the volume $V_{t}$ of an arbitrary chain 
bounding $f_{t}$ in $M_{t}$. Then,
$$\dot V_{0} = {\textstyle {1\over 2}} l_{0}\bigl 
( \dot b_{0}\bigr )  ,  \tag  3$$
where $b_{t}\in {\Cal H}\left ( \lambda ;{\Bbb 
R}/2\pi {\Bbb Z}\right )$ is the bending cocycle 
of the pleated surface $f_{t}$, and where the 
right hand term denotes one half of the length of 
$\dot b_{0}\in {\Cal H}\left ( \lambda ;{\Bbb 
R}\right )$ with respect to the hyperbolic metric 
$m_{0}$ on $S$ defined by pull back under $f_{0}$ 
of the hyperbolic metric of $M_{0}$.  
\endproclaim    
\demo {Proof of Theorem~3 when there are no 
cusps}   As traditional in 3--dimensional 
hyperbolic geometry, the presence of cusps 
introduces some local technicalities which are 
not difficult, but tend to dilute attention away 
from the main points of the proof. For this 
reason, we will first restrict ourselves to the 
case where the hyperbolic manifolds $M_{t}$ have 
no cusps, and we will later explain how to extend 
the proof in the presence of cusps.
\nobreak\par\nobreak 	Consequently, assume that 
the surface $S$ is compact and that the manifolds 
$M_{t}$ have no cusps. 
\nobreak\par\nobreak 	For every $t$, let $m_{t}$ 
be the hyperbolic metric on $S$ obtained by 
pulling back the metric of $M_{t}$ under $f_{t}$, 
and let $\lambda _{t}$ denote the 
$m_{t}$--geodesic lamination of $S$ corresponding 
to the geodesic lamination $\lambda $. By 
hypothesis, $\lambda $ is a maximal geodesic 
lamination and the complement $S-\lambda _{t}$ 
consists of ideal triangles. 
\nobreak\par\nobreak 	We can cover $\lambda _{0}$ 
by a family of rectangles $R_{1}^{(0)}$, 
$R_{2}^{(0)}$, \dots , $R_{m}^{(0)}$ with 
$m_{0}$--geodesic sides, with disjoint interiors, 
and such that the components of $\lambda _{0}\cap 
R_{i}^{(0)}$ are all parallel to (and disjoint 
from) two opposite sides of the rectangle 
$R_{i}^{(0)}$, for each $i$. These rectangles 
more or less form a train track carrying $\lambda 
_{0}$. If we collapse each rectangle 
$R_{i}^{(0)}$ to the edge that is parallel to the 
components of $\lambda _{0}\cap R_{i}^{(0)}$, we 
obtain a graph embedded in $S$. Extend this graph 
to a triangulation ${\Cal T}$ of $S$ and choose a 
map $g_{0}:S\rightarrow M_{0}$ which is homotopic 
to $f_{0}$, is polyhedral  with respect to ${\Cal 
T}$, and sends to a geodesic arc the image of 
each rectangle $R_{j}^{(0)}$ in the 1--skeleton 
of ${\Cal T}$. 
\nobreak\par\nobreak 	For $t$ small, we similarly 
construct rectangles $R_{1}^{(t)}$, 
$R_{2}^{(t)}$, \dots , $R_{m}^{(t)}$ with 
$m_{t}$--geodesic sides, with disjoint interiors, 
and such that the components of each $\lambda 
_{t}\cap R_{i}^{(t)}$ are parallel to one side of 
$R_{i}^{(t)}$. In addition, we require these 
$R_{i}^{(t)}$ to vary differentiably with $t$. In 
particular, collapsing the $R_{i}^{(t)}$ defines 
the same graph embedded in $S$, up to isotopy. 
Choose a map  $g_{t}:S\rightarrow M_{t}$ that is 
homotopic to $f_{t}$, is polyhedral with respect 
to the triangulation ${\Cal T}$, sends the image 
of each rectangle $R_{j}^{(t)}$ under the 
collapsing process to a geodesic arc, and varies 
differentiably with $t$. 
\nobreak\par\nobreak 	The Schl{\"a}fli formula of 
Corollary~2 determines the variation of the 
volume enclosed by the polyhedral map $g_{t}$. To 
compute the variation of the volume enclosed by 
$f_{t}$, it is therefore sufficient to analyze 
the volume of a homotopy between $f_{t}$ and 
$g_{t}$. In the definition of the triangulation 
${\Cal T}$, we implicitly used a map 
$h_{t}:S\rightarrow S$ homotopic to the identity 
and collapsing each rectangle $R_{i}^{(t)}$ to an 
arc contained in the 1-skeleton of ${\Cal T}$. 
Since there is a volume 0 homotopy between 
$g_{t}$ and $g_{t}\circ h_{t}$, it suffices to 
determine the volume of a homotopy $H_{t}:S\times 
\left \lbrack 0,1\right \rbrack \rightarrow 
M_{t}$ such that $H_{t \vert S\times \left 
\lbrace 1\right \rbrace }=f_{t}$ and $H_{t \vert 
S\times \left \lbrace 0\right \rbrace 
}=g_{t}\circ h_{t}$. 
\nobreak\par\nobreak 	Let us focus attention on a 
rectangle $R_{i}^{(t)}$. We `straighten' the 
homotopy $H_{t}$ on $R_{i}^{(t)}\times \left 
\lbrack 0,1\right \rbrack $ in the following way. 
Identify $R_{i}^{(t)}$ to a standard rectangle 
$\left \lbrack a,b\right \rbrack \times \left 
\lbrack c,d\right \rbrack $ by an 
orientation-preserving homeomorphism such that 
each component of $\lambda _{t}\cap R_{i}^{(t)}$ 
corresponds to an arc $\left \lbrace x\right 
\rbrace \times \left \lbrack c,d\right \rbrack $. 
 Cut each rectangle $\left \lbrace x\right 
\rbrace \times \left \lbrack c,d\right \rbrack 
\times \left \lbrack 0,1\right \rbrack $ in 
$R_{i}^{(t)}\times \left \lbrack 0,1\right 
\rbrack $ into two triangles along the diagonal 
line joining $\left ( x,c,1\right )$ to $\left ( 
x,d,0\right )$, foliate the upper triangle by 
line segments originating from $\left ( 
x,d,0\right )$, and foliate the lower rectangle 
by line segments originating from $\left ( 
x,c,0\right )$. This decomposes 
$R_{i}^{(t)}\times \left \lbrack 0,1\right 
\rbrack \cong \left \lbrack a,b\right \rbrack 
\times \left \lbrack c,d\right \rbrack \times 
\left \lbrack 0,1\right \rbrack $ into a family 
of arcs; see Figure~1. We can now deform the 
restriction of $H_{t}$ to $R_{i}^{(t)}\times 
\left \lbrack 0,1\right \rbrack $ so that it 
sends each of these arcs to a geodesic arc.

\midinsert
\TrimBottom{2pc}
\centerline{\BoxedEPSF{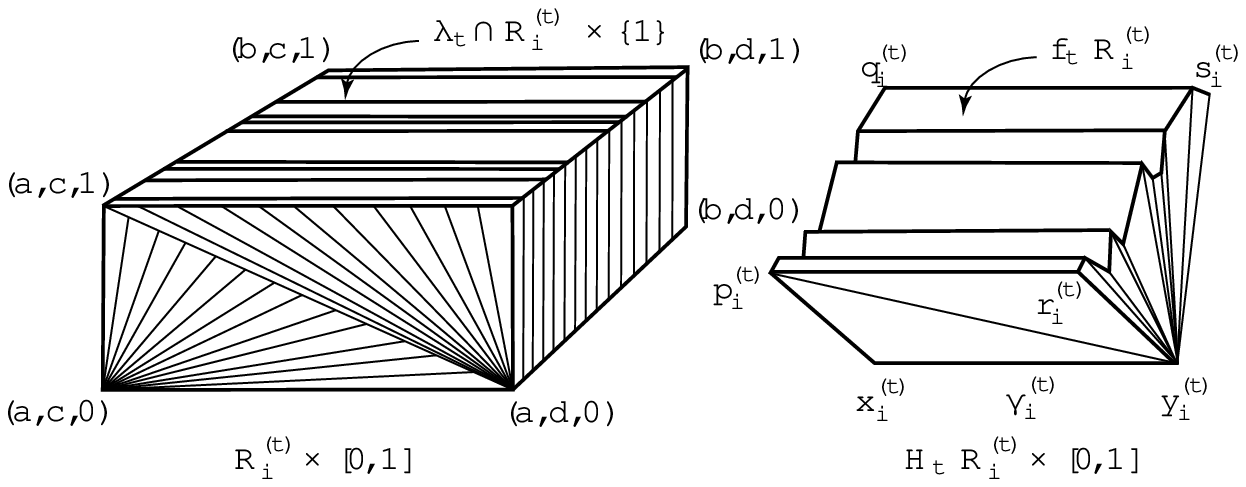}}
\botcaption {Figure~1} \endcaption
\endinsert

\nobreak\par\nobreak 	By construction, 
$H_{t}\bigl ( R_{i}^{(t)}\times \left \lbrace 
1\right \rbrace \bigr ) $ is equal to $f_{t}\bigl 
( R_{i}^{(t)}\bigr ) $, and $H_{t}\bigl ( 
R_{i}^{(t)}\times \left \lbrace 0\right \rbrace 
\bigr ) =g_{t}\circ h_{t}\bigl ( R_{i}^{(t)}\bigr 
) $ is a geodesic arc with end points 
$x_{i}^{(t)}=H_{t}\left ( \left \lbrack a,b\right 
\rbrack \times \left \lbrace c\right \rbrace 
\times \left \lbrace 0\right \rbrace \right )$ 
and $y_{i}^{(t)}=H_{t}\left ( \left \lbrack 
a,b\right \rbrack \times \left \lbrace d\right 
\rbrace \times \left \lbrace 0\right \rbrace 
\right )$. The face $\left \lbrack a,b\right 
\rbrack \times \left \lbrace c\right \rbrace 
\times \left \lbrack 0,1\right \rbrack $ of 
$R_{i}^{(t)}\times \left \lbrack 0,1\right 
\rbrack $ is sent by $H_{t}$ to a `pleated fan' 
which is the joint of the arc $f_{t}\left ( \left 
\lbrack a,b\right \rbrack \times \left \lbrace 
c\right \rbrace \right )$ and of the point 
$x_{i}^{(t)}$, namely $H_{t}\left ( \left \lbrack 
a,b\right \rbrack \times \left \lbrace c\right 
\rbrace \times \left \lbrack 0,1\right \rbrack 
\right )$ consists of all geodesic arcs that join 
$f_{t}\left ( \left \lbrack a,b\right \rbrack 
\times \left \lbrace c\right \rbrace \right )$ to 
$x_{i}^{(t)}$ and are in the appropriate homotopy 
class. Similarly,  $\left \lbrack a,b\right 
\rbrack \times \left \lbrace d\right \rbrace 
\times \left \lbrack 0,1\right \rbrack $ is sent 
to the joint of $f_{t}\left ( \left \lbrack 
a,b\right \rbrack \times \left \lbrace d\right 
\rbrace \right )$ and $y_{i}^{(t)}$. The 
remaining two faces $\left \lbrace a\right 
\rbrace \times \left \lbrack c,d\right \rbrack 
\times \left \lbrack 0,1\right \rbrack $ and 
$\left \lbrace b\right \rbrace \times \left 
\lbrack c,d\right \rbrack \times \left \lbrack 
0,1\right \rbrack $ of $R_{i}^{(t)}\times \left 
\lbrack 0,1\right \rbrack $ are each sent to the 
union of two totally geodesic triangles.
\nobreak\par\nobreak 	In particular, this 
analysis of the restriction of $H_{t}$ to the 
faces $\left \lbrack a,b\right \rbrack \times 
\left \lbrace c\right \rbrace \times \left 
\lbrack 0,1\right \rbrack $ and $\left \lbrack 
a,b\right \rbrack \times \left \lbrace d\right 
\rbrace \times \left \lbrack 0,1\right \rbrack $ 
shows that we can globally deform $H_{t}$ so that 
it is of the above type on each rectangle 
$R_{i}^{(t)}$. 
\nobreak\par\nobreak 	To evaluate the volume of 
the restriction of $H_{t}$ to $R_{i}^{(t)}\times 
\left \lbrack 0,1\right \rbrack $, we decompose 
it into pieces. For each component $R$ of 
$R_{i}^{(t)}-\lambda _{t}$, $H_{t}\left ( R\times 
\left \lbrack 0,1\right \rbrack \right )$ is the 
union of a pyramid with square basis, namely the 
joint of $f_{t}\left ( R\right )$ and 
$y_{i}^{(t)}$, and of the tetrahedron formed by 
the joint of the two geodesic arcs $f_{t}\left ( 
R\cap \left \lbrack a,b\right \rbrack \times 
\left \lbrace c\right \rbrace \right )$ and 
$\gamma _{i}^{(t)}$. For each component $k$ of 
$\lambda _{t}\cap R_{i}^{(t)}$, $H_{t}\left ( 
k\times \left \lbrack 0,1\right \rbrack \right )$ 
is the union of the two totally geodesic 
triangles which are, respectively, the joint of 
$f_{t}\left ( k\right )$ and $y_{i}^{(t)}$ and 
the joint of $f_{t}\left ( k\cap \left ( \left 
\lbrack a,b\right \rbrack \times \left \lbrace 
c\right \rbrace \right )\right )$ and $\gamma 
_{i}^{(t)}$. For the metric $m_{t}$, $\lambda 
_{t}$ has 2--dimensional Lebesgue measure 0 in 
$S$, and $\lambda _{t}\cap \left \lbrack 
a,b\right \rbrack \times \left \lbrace c\right 
\rbrace $ has 1--dimensional Lebesgue measure 0 
in the transverse arc $\left \lbrack a,b\right 
\rbrack \times \left \lbrace c\right \rbrace $ 
\cite{CaB}\cite{PeH}\cite{BiS}. It follows that 
$H_{t}\bigl ( \bigl ( \lambda _{t}\cap 
R_{i}^{(t)}\bigr ) \times \left \lbrack 0,1\right 
\rbrack \bigr ) $ has 3--dimensional Lebesgue 
measure 0. Therefore, we can focus only on the 
contribution of the components of 
$R_{i}^{(t)}-\lambda _{t}$. 
\nobreak\par\nobreak 	We can now sketch the proof 
of Theorem~3. Let $R^{(t)}\subset S$ denote the 
union of the rectangles $R_{i}^{(t)}$. By 
construction, $H_{t}\left ( \left ( 
S-R^{(t)}\right )\times \left \lbrack 0,1\right 
\rbrack \right )$ is bounded by a polyhedral 
surface, and the variation of its volume is given 
by Corollary~2. We observed that the volume of 
$H_{t}\left ( R^{(t)}\times \left \lbrack 
0,1\right \rbrack \right )$ is equal to the sum 
of the volumes of certain pyramids and 
tetrahedra. We can therefore expect that the 
variation of the volume of $H_{t}\left ( 
R^{(t)}\times \left \lbrack 0,1\right \rbrack 
\right )$ is equal to the sum of the variations 
of the volumes of these pyramids and tetrahedra, 
as given by Corollary~2. Altogether, this 
expresses the variation of the volume of $H_{t}$ 
as the sum of lengths of edges  multiplied by the 
variation of dihedral angles at these edges. As 
in the proof of Corollary~2, the contributions of 
the internal edges cancel out, as well as the 
contribution of the edges that are contained in 
the sides of the rectangles $R_{i}^{(t)}$. This 
leaves the contributions of the edges of the 
polyhedral map $g_{t}$, which will cancel out 
with the variation of the volume enclosed by 
$g_{t}$, and the contribution of the edges 
contained in $f_{t}\left ( \lambda _{t}\right )$, 
which can be re-interpreted as the length of the 
variation of the bending cocycle of the pleated 
surface $f_{t}$.
\nobreak\par\nobreak 	These ideas easily lead to 
a formal proof of Theorem~3, but numerous points 
need to be justified. First of all, to be able to 
apply Corollary~2, we need to know that the 
shapes of the pyramids and tetrahedra of the 
decomposition are non-degenerate and vary 
differentiably with $t$. Then, because there are 
(in general) infinitely many such pyramids and 
tetrahedra, we have to show that the infinite 
sums involved do converge. Because the internal 
edges are not locally finite, the proof that 
their contributions cancel out is not as simple 
as in the proof of Corollary~2. Finally, we have 
to identify the contribution of the edges 
contained in $f_{t}\left ( \lambda _{t}\right )$ 
to the length of the variation of the bending 
cocycle of the pleated surface $f_{t}$.

\proclaim {Lemma 4 }   Given a component $R_{0}$ 
of $R_{i}^{(0)}-\lambda _{0}$, let $R_{t}$ be the 
corresponding component of $R_{i}^{(t)}-\lambda 
_{t}$.  Then, the vertices of the rectangle 
$f_{t}\left ( R_{t}\right )$ vary differentiably 
with $t$ in $M_{t}$.  \endproclaim    
\demo {Proof of Lemma~4}   Recall that we are 
given a family of representations $\rho 
_{t}:\Gamma \rightarrow \roman {Isom}^{+}\left ( 
{\Bbb H}^{3}\right )$ depending differentiably on 
$t$ and of isometries $M_{t}\cong {\Bbb 
H}^{3}/\rho _{t}\left ( \Gamma \right )$. We want 
to show that $f_{t}\left ( R_{t}\right )\subset 
M_{t}\cong {\Bbb H}^{3}/\rho _{t}\left ( \Gamma 
\right )$ lifts to a rectangle depending 
differentiably on $t$ in ${\Bbb H}^{3}$. 
\nobreak\par\nobreak 	In \cite{Bo3, {\S }10}, we 
showed that the restriction of $f_{t}$ to each 
component of $S-\lambda _{t}$ depends 
differentiably on $t$. Namely, if we lift $f_{t}$ 
to a pleated surface $\widetilde  
f_{t}:\widetilde  S\rightarrow {\Bbb H}^{3}$ 
defined on the universal covering $\widetilde  
S$, then for every component $P$ of the 
complement of the preimage $\widetilde  \lambda $ 
of $\lambda $ in $\widetilde  S$ and if $P_{t}$ 
denotes the corresponding component of 
$\widetilde  S-\widetilde  \lambda _{t}$, the 
ideal triangle $\widetilde  f_{t}\left ( 
P_{t}\right )$ in ${\Bbb H}^{3}$ depends 
differentiably on $t$. (Strictly speaking, we 
proved this property only if we replace the 
isometric identification $M_{t}\cong {\Bbb 
H}^{3}/\rho _{t}\left ( \Gamma \right )$ by 
another identification $M_{t}\cong {\Bbb 
H}^{3}/\rho '_{t}\left ( \Gamma \right )$ where 
there exists $A_{t}\in \roman {Isom}^{+}\left ( 
{\Bbb H}^{3}\right )$ such that $\rho _{t}'\left 
( \gamma \right )=A_{t}^{-1}\rho _{t}\left ( 
\gamma \right )A_{t}$ for every $\gamma \in 
\Gamma $. For every $\gamma \in \pi _{1}\left ( 
S\right )\subset \Gamma $, the fact that 
$\widetilde  f_{t}\left ( P_{t}\right )$ and 
$\widetilde  f_{t}\left ( \gamma P_{t}\right )$ 
depend differentiably on $t$ for an arbitrary 
component $P$ of $\widetilde  S-\widetilde  
\lambda $ show that $\rho '_{t}\left ( \gamma 
\right )$ depends differentiably on $t$. Looking 
at the fixed points of the isometry groups $\rho 
_{t}\left ( \pi _{1}\left ( S\right )\right )$ 
and $\rho '_{t}\left ( \pi _{1}\left ( S\right 
)\right )$, we conclude that $A_{t}$ depends 
differentiably on $t$. Using the identification 
${\Bbb H}^{3}/\rho _{t}\left ( \Gamma \right 
)\cong {\Bbb H}^{3}/\rho _{t}'\left ( \Gamma 
\right )$ induced by $A_{t}$, we can therefore 
assume that $\rho _{t}'=\rho _{t}$.)
\nobreak\par\nobreak 	\cite{Bo3} also shows that 
the pull back metric $m_{t}$ on $S$ depends 
differentiably on $t$. Again, \cite{Bo3, {\S }5} 
provides a representation $\sigma _{t}:\pi 
_{1}\left ( S\right )\rightarrow \roman 
{Isom}^{+}\left ( {\Bbb H}^{2}\right )$ depending 
differentiably on $t$ and an isometric 
identification $\varphi _{t}:\left ( 
S,m_{t}\right )\rightarrow {\Bbb H}^{2}/\sigma 
_{t}\left ( \pi _{1}\left ( S\right )\right )$ 
such that, for every component $P$ of $\widetilde 
 S-\widetilde  \lambda $, the image $\widetilde  
\varphi _{t}\left ( P_{t}\right )$ of the 
corresponding component $P_{t}$ of $\widetilde  
S-\widetilde  \lambda _{t}$ under a lift 
$\widetilde  \varphi _{t}:\widetilde  
S\rightarrow {\Bbb H}^{2}$ is an ideal triangle 
which varies differentiably with $t$. We assumed 
that the vertices of the rectangles $R_{i}^{(t)}$ 
depend differentiably on $t$ for the metric 
$m_{t}$, namely that, if we lift $R_{i}^{(t)}$ to 
$\widetilde  R_{i}^{(t)}\subset \widetilde  S$, 
the vertices of $\widetilde  \varphi _{t}\bigl ( 
\widetilde  R_{i}^{(t)}\bigr ) \subset {\Bbb 
H}^{2}$ depend differentiably on $t$. If $P_{t}$ 
is the component of $\widetilde  S-\widetilde  
\lambda _{t}$ that contains the lift $\widetilde  
R_{t}\subset \widetilde  R_{i}^{(t)}$ of $R_{t}$, 
and if we $m_{t}$--isometrically identify $P_{t}$ 
to a fixed ideal triangle, it follows that 
$\widetilde  R_{t}=\widetilde  \varphi 
_{t}^{-1}\bigl ( \widetilde  \varphi _{t}\bigl ( 
\widetilde  R_{i}^{(t)}\bigr ) \cap \widetilde  
\varphi _{t}\left ( P_{t}\right )\bigr ) $ 
depends differentiably on $t$ in this fixed ideal 
triangle.
\nobreak\par\nobreak 	Because $\widetilde  
f_{t}\left ( P_{t}\right )$ depends 
differentiably on $t$, we conclude that 
$\widetilde  f_{t}\bigl ( \widetilde  R_{t}\bigr 
) $ depends differentiably on $t$ in ${\Bbb 
H}^{3}$, and therefore that $f_{t}\left ( 
R_{t}\right )$ depends differentiably on $t$ in 
$M_{t}$.  \qed\enddemo

\nobreak\par\nobreak 	Lemma~4 shows that each of 
the pyramids and tetrahedra of the decomposition 
of $H_{t}\bigl ( R_{i}^{(t)}\times \left \lbrack 
0,1\right \rbrack \bigr ) $ varies differentiably 
with $t$. 
\nobreak\par\nobreak 	We need to precise Lemma~4, 
using the estimates of \cite{Bo3}. When $R_{t}$ 
does not contain one of the sides of 
$R_{i}^{(t)}$, the two leaves of $\lambda _{t}$ 
that it touches follow each other for a while, 
crossing the same rectangles $R_{j}^{(t)}$. 
However, in some direction, they must diverge and 
cross different $R_{j}^{(t)}$ after a while, 
since they would otherwise stay within bounded 
distance of each other and therefore be equal. 
Let the {\bi divergence radius \/} $r\left ( 
R_{t}\right )\geqslant 1$ be the number of 
$R_{j}^{(t)}$ which they cross in common before 
diverging. By convention, $r\left ( R_{t}\right 
)=1$ when $R_{t}$ contains one of the sides of 
$R_{i}^{(t)}$.

\proclaim {Lemma~5 }   With the data of Lemma~4, 
the derivative $\dot a_{t}$ of each vertex 
$a_{t}$ of $f_{t}\left ( R_{t}\right )$ with 
respect to $t$ is an $O\left ( r\left ( 
R_{t}\right )\right )=O\left ( r\left ( 
R_{0}\right )\right )$. In addition, for any two 
vertices $a_{t}$, $b_{t}$, the distance between 
the vectors $\dot a_{t}$ and $\dot b_{t}$ is an 
$O\left ( d\left ( a_{t},b_{t}\right )r\left ( 
R_{t}\right )\right )$.     \endproclaim    
\demo {Proof}   To give a sense to this 
statement, we need to choose a lift of 
$f_{t}\left ( R_{t}\right )\subset M_{t}\cong 
{\Bbb H}^{3}/\rho _{t}\left ( \Gamma \right )$ to 
${\Bbb H}^{3}$, as in Lemma~4. For this, we first 
choose a lift of $f_{t}\bigl ( R_{i}^{(t)}\bigr ) 
$ to ${\Bbb H}^{3}$, and then restrict it to a 
lift of each  $f_{t}\left ( R_{t}\right )$. The 
statement of Lemma~5 implicitly assumes that the 
lifts of the $f_{t}\left ( R_{t}\right )$ are 
chosen in such a consistent way. The constants 
hidden in the symbols $O\left ( \enskip \right )$ 
will then depend on the choice of the lift of 
$f_{t}\bigl ( R_{i}^{(t)}\bigr ) $, but not on 
the components $R_{t}$.
\nobreak\par\nobreak 	Let us use the notation of 
the proof of Lemma~4. In \cite{Bo3, {\S }10}, we 
give explicit formulas expressing the restriction 
of $\widetilde  f_{t+h}$ to each component of 
$\widetilde  S-\widetilde  \lambda _{t+h}$ as a 
limit of rotation-translations along geodesics of 
$\widetilde  f_{t}\bigl ( \widetilde  \lambda 
_{t}\bigr ) $. In addition, we show that the 
convergence is holomorphic, so that we can 
differentiate in the limit. Differentiating with 
respect to $h$ and applying to the derivative the 
estimates of \cite{Bo3, {\S }5} and \cite{Bo3, 
Lemma~6}, it easily follows that the derivative 
of the ideal triangle $\widetilde  f_{t}\left ( 
P_{t}\right )$ is an $O\left ( r\left ( 
R_{t}\right )\right )$. Similarly, the derivative 
of $\widetilde  \varphi _{t}\left ( P_{t}\right 
)$ is also an $O\left ( r\left ( R_{t}\right 
)\right )$. The first statement of Lemma~5 
immediately follows.
\nobreak\par\nobreak 	For the second statement, 
we can restrict attention to the case where 
$a_{t}$ and $b_{t}$ are on the same side of 
$f_{t}\bigl ( R_{i}^{(t)}\bigr ) $; indeed, 
$d\left ( a_{t},b_{t}\right )$ is otherwise 
bounded away from 0, and the result is trivial. 
Then, $a_{t}$ and $b_{t}$ are the images under 
$\widetilde  f_{t}\circ \bigl ( \widetilde  
\varphi _{t}{}_{ \vert P_{t}}\bigr ) ^{-1}$ of 
the intersection points of $\widetilde  \varphi 
_{t}\left ( \partial P_{t}\right )$ with an arc 
$k_{t}$ in ${\Bbb H}^{2}$ which varies 
differentiably with $t$ (and corresponds to a 
side of $R_{i}^{(t)}$). Because $\widetilde  
f_{t}{}_{ \vert P_{t}}$ and $\widetilde  \varphi 
_{t}{}_{ \vert P_{t}}$ are isometries, their 
differentials are uniformly Lipschitz, and the 
$\roman {C}^{1}$--norm of their derivatives with 
respect to $t$ is an $O\left ( r\left ( 
R_{t}\right )\right )$ by the previous estimate. 
The second statement of Lemma~ 5 then follows 
from the chain rule.  \qed\enddemo

\nobreak\par\nobreak 	To apply the Schl{\"a}fli 
formula to these pyramids and tetrahedra, we need 
to make sure that their faces are not collapsed 
to arcs or points (otherwise, dihedral angles do 
not make sense). This means that for every $t$ 
the points $x_{i}^{(t)}$ and $y_{i}^{(t)}$ avoid 
a countable union of geodesic arcs of bounded 
length in $M_{t}$. By Lemma~4, these geodesic 
arcs depend differentiably on $t$. As a 
consequence, their lifts sweep a domain of 
Lebesgue measure 0 in ${\Bbb H}^{3}$. We can 
therefore choose the polyhedral maps $g_{t}$ 
generic enough so that the pyramids and 
tetrahedra are never degenerate. Actually, we 
only need this property for the sake of the 
exposition, since the estimate of Lemma~6 below 
would enable us to deal with degenerate pyramids 
and tetrahedra as well. 
\nobreak\par\nobreak 	Also, the image $f_{t}\left 
( \lambda _{t}\right )$ of the pleating locus has 
Hausdorff dimension 1 \cite{BiS}, and varies 
continuously with $t$; see \cite{Thu, {\S 
}8}\cite{CEG, {\S }5} or \cite{Bo3}. Therefore, 
we can arrange that the points $x_{i}^{(t)}$ and 
$y_{i}^{(t)}$ stay at distance bounded away from 
0 from $f_{t}\left ( \lambda _{t}\right )$.  
\nobreak\par\nobreak 	Let $R_{t}$ be a component 
of $R_{i}^{(t)}-\lambda _{t}$, and let $p_{t}$, 
$q_{t}$, $r_{t}$, $s_{t}$ be the vertices of the 
totally geodesic rectangle $f_{t}\left ( 
R_{t}\right )$ where, for the identification 
$R_{i}^{(t)}\cong \left \lbrack a,b\right \rbrack 
\times \left \lbrack c,d\right \rbrack $, the 
points $p_{t}$, $q_{t}$ occur in this order in 
$f_{t}\left ( \left \lbrack a,b\right \rbrack 
\times \left \lbrace c\right \rbrace \right )$ 
and the points $r_{t}$, $s_{t}$ occur in this 
order in $f_{t}\left ( \left \lbrack a,b\right 
\rbrack \times \left \lbrace d\right \rbrace 
\right )$. Then, the pyramid $P_{t}$ associated 
to $R_{t}$ has vertices $p_{t}$, $q_{t}$, 
$r_{t}$, $s_{t}$ and $y_{i}^{(t)}$. See Figure~2.

\midinsert
\centerline{\BoxedEPSF{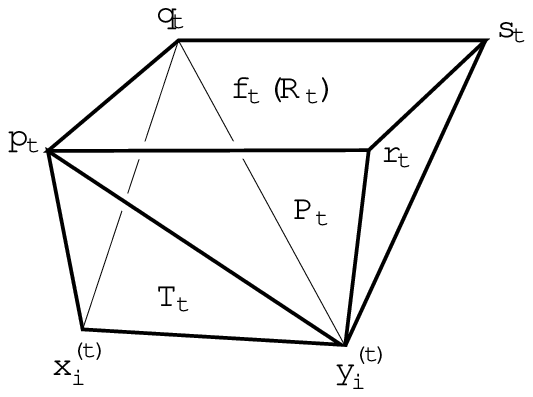}}
\botcaption {Figure~2} \endcaption
\endinsert

\nobreak\par\nobreak 	The Schl{\"a}fli formula 
gives that $d\roman {vol}\left ( P_{t}\right ) 
/dt$ is the sum of the terms
$${\textstyle {1\over 2}} l\left ( 
p_{t}q_{t}\right ) {d\over dt} \theta 
_{P_{t}}\bigl ( p_{t}q_{t}\bigr ) , \tag  4$$
$${\textstyle {1\over 2}} l\left ( 
r_{t}s_{t}\right ) {d\over dt} \theta 
_{P_{t}}\bigl ( r_{t}s_{t}\bigr )  , \tag  5$$
$${\textstyle {1\over 2}} \bigl ( l\bigl ( 
p_{t}y_{i}^{(t)}\bigr )  {d\over dt} \theta 
_{P_{t}}\bigl ( p_{t}y_{i}^{(t)}\bigr )  + l\bigl 
( q_{t}y_{i}^{(t)}\bigr )  {d\over dt} \theta 
_{P_{t}}\bigl ( q_{t}y_{i}^{(t)}\bigr )  \bigr ) 
, \tag  6$$
$${\textstyle {1\over 2}} \bigl ( l\bigl ( 
r_{t}y_{i}^{(t)}\bigr )  {d\over dt} \theta 
_{P_{t}}\bigl ( r_{t}y_{i}^{(t)}\bigr )  + l\bigl 
( s_{t}y_{i}^{(t)}\bigr )  {d\over dt} \theta 
_{P_{t}}\bigl ( s_{t}y_{i}^{(t)}\bigr )  \bigr ) 
, \tag 7 $$
and
$${\textstyle {1\over 2}} \bigl ( l\left ( 
p_{t}r_{t}\right ) {d\over dt} \theta 
_{P_{t}}\bigl ( p_{t}r_{t}\bigr )  + l\left ( 
s_{t}q_{t}\right ) {d\over dt} \theta 
_{P_{t}}\bigl ( s_{t}q_{t}\bigr )  \bigr )   \tag 
 8$$
where $l\left ( \enskip \right )$ denotes the 
length of the edge indicated, and where $\theta 
_{P_{t}}\left ( \enskip \right )$ is the external 
dihedral angle of the boundary $\partial P_{t}$ 
at the edge indicated. For this, we orient 
$\partial P_{t}$ so that the restriction $f_{t 
\vert P_{t}}:R_{t}\rightarrow f_{t}\left ( 
R_{t}\right )\subset \partial P_{t}$ is 
orientation -preserving and we orient $P_{t}$ 
accordingly (so that $\roman {vol}\left ( 
P_{t}\right )$ may be negative).  
\nobreak\par\nobreak 	The grouping of the terms 
is here important, because it will guarantee the 
convergence of the series when we sum over all 
components $R_{t}$ of $R_{i}^{(t)}-\lambda _{t}$. 
It is not hard to see that this sum will not 
converge if we do not use this grouping.
\nobreak\par\nobreak 	We first sum the terms of 
type \rom{(4)}. By the second part of Lemma~5, 
$d\theta _{P_{t}}\bigl ( p_{t}q_{t}\bigr ) /dt = 
O\left ( r\left ( R_{0}\right )+1\right )$. Also, 
there is a constant $A>0$, depending only on the 
length of the components of $\lambda _{t}\cap 
R_{i}^{(t)}$ (and therefore uniform in $t$), such 
that the leaves of $\lambda _{t}$ passing through 
$p_{t}$ and $q_{t}$ stay at uniformly bounded 
distance from each other over a length of at 
least $Ar\left ( R_{t}\right )=Ar\left ( 
R_{0}\right )$; since the metric $m_{t}$ is 
hyperbolic, it follows that $l\left ( 
p_{t}q_{t}\right )=O\left ( e^{-Ar\left ( 
R_{0}\right )}\right )$. Finally, for every 
$r\geqslant 1$, the number of components $R_{0}$ 
of $R_{i}^{(0)}-\lambda _{0}$ such that $r\left ( 
R_{0}\right )=r$ is uniformly bounded by the 
number of spikes of $S-\lambda _{0}$. It follows 
that, as we sum over all components $R_{t}$ of 
$R_{i}^{(t)}-\lambda _{t}$, the series
$${\textstyle {1\over 2}} \sum _{R_{t}} l\left ( 
p_{t}q_{t}\right ) {d\over dt} \theta 
_{P_{t}}\bigl ( p_{t}q_{t}\bigr )  \tag  9$$
is convergent. 
\nobreak\par\nobreak 	The same arguments show the 
convergence of the terms of type \rom{(5)}, 
namely
$${\textstyle {1\over 2}} \sum _{R_{t}} l\left ( 
r_{t}s_{t}\right ) {d\over dt} \theta 
_{P_{t}}\bigl ( r_{t}s_{t}\bigr )  . \tag  10$$
\nobreak\par\nobreak 	To show the convergence of 
the sum of the terms of type \rom{(6)}, we first 
estimate the quantity 
$$d \theta _{P_{t}}\bigl ( p_{t}y_{i}^{(t)}\bigr 
)  /dt+ d \theta _{P_{t}}\bigl ( 
q_{t}y_{i}^{(t)}\bigr ) /dt .$$ 
For this, choose an isometric embedding $\varphi 
_{t}:P_{t}\rightarrow {\Bbb H}^{3}$ such that 
$\varphi _{t}\left ( f_{t}\left ( R_{t}\right 
)\right )$ is contained in ${\Bbb H}^{2}\subset 
{\Bbb H}^{3}$, $\varphi _{t}\left ( 
p_{t}r_{t}\right )$ is contained in a fixed 
geodesic $g$ of ${\Bbb H}^{2}$, and  $\varphi 
_{t}\left ( p_{t}\right )$ is a fixed point $p\in 
g$. Since the two leaves of $\lambda _{t}$ 
touching $R_{t}$ are asymptotic, the geodesic of 
${\Bbb H}^{2}$ that contains $\varphi _{t}\left ( 
q_{t}\right )$ and $\varphi _{t}\left ( 
s_{t}\right )$ has a common end point $x$ with 
$g$ at infinity. We can now consider the 
tetrahedron $T$ with vertices $p=\varphi 
_{t}\left ( p_{t}\right )$, $q=\varphi _{t}\left 
( q_{t}\right )$, $y=\varphi _{t}\bigl ( 
y_{i}^{(t)}\bigr ) $ and $x$; see Figure~3.

\midinsert
\centerline{\BoxedEPSF{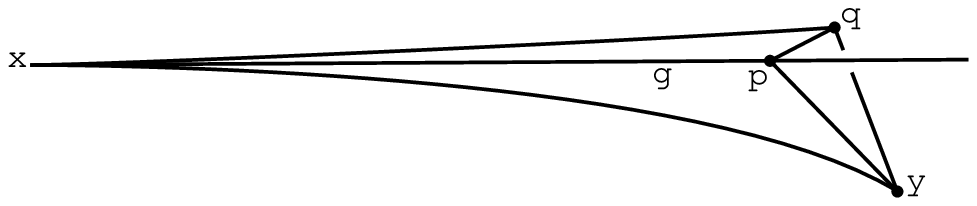}}
\botcaption {Figure~3} \endcaption
\endinsert

\proclaim {Lemma~6 }   Let $g$ be a geodesic of 
${\Bbb H}^{2}\subset {\Bbb H}^{3}$, let $p$ be a 
fixed point of $g$, and let $x$ be one of the end 
points of $g$ on the circle at infinity $\partial 
_{\infty }{\Bbb H}^{2}$. Given two constants 
$A>0$ and $B>0$, consider two points $q\in {\Bbb 
H}^{2}$ and $y\in {\Bbb H}^{3}$ such that the 
distances from $y$ to $p$ and from $p$ to $q$ are 
at most $A$, and such that the distance from $y$ 
to $g$ and to the geodesic of ${\Bbb H}^{2}$ 
containing $q$ and $x$ is at least $B$. Finally, 
in the tetrahedron $T$ of vertices $p$, $q$, $x$, 
$y$, let $\Theta \left ( q,y\right )$ denote the 
sum of the internal dihedral angles of $T$ along 
the edges $py$ and $qy$. Then, $\Theta \left ( 
q,y\right )$ is a differentiable function of $q$ 
and $y$. In addition, if $q$ varies with velocity 
$\dot q$ and $y$ varies with velocity $\dot y$, 
the derivative $\dot \Theta \left ( q,y\right )$ 
of $\Theta \left ( q,y\right )$ is an $O\left ( 
\left \Vert \dot q\right \Vert +d\left ( 
p,q\right )\left \Vert \dot y\right \Vert \right 
)$, where the constant hidden in the symbol 
$O\left ( \enskip \right )$ depends only on the 
constants $A$ and $B$.  \endproclaim

\demo {Proof}   Note that our definition of 
$\Theta \left ( q,y\right )$ does not make sense 
when the tetrahedron $T$ is degenerate, namely 
when the three points $p$, $q$ and $y$ are on the 
same geodesic (since $B>0$, the triangles $xpy$ 
and $xqy$ cannot be degenerate). We first extend 
it to this case, by using a different point of 
view.
\nobreak\par\nobreak 	At the point $y$, consider 
the unit vectors  $v_{p}$, $v_{q}$, $v_{x}$ 
pointing in the direction of $p$, $q$, $x$, 
respectively. These vectors draw a triangle 
$v_{x}v_{p}v_{q}$ on the visual sphere at $y$. 
Then, $\Theta \left ( q,y\right )$ is the sum of 
the angles of this spherical triangle at $v_{p}$ 
and $v_{q}$, when these angles make sense. By the 
Gauss formula, $\Theta \left ( q,y\right )$ is 
therefore equal to $\pi $ plus the area of the 
spherical triangle $v_{x}v_{p}v_{q}$ minus the 
angle of the triangle at $v_{x}$. Since $y$ stays 
away from $g$ and the geodesic containing $x$ and 
$q$, the edges $v_{x}v_{p}$ and $v_{x}v_{q}$ of 
the triangle are never reduced to a point. This 
formula for $\Theta \left ( q,y\right )$ 
consequently makes sense for every $q$, $y$ 
satisfying the conditions of the Lemma, and shows 
that $\Theta \left ( q,y\right )$ is an 
infinitely differentiable function of $q$ and 
$y$.
\nobreak\par\nobreak 	The positions allowed by 
the conditions of Lemma~6 for the point $\left ( 
q,y\right )$ form a compact subset of ${\Bbb 
H}^{2}\times {\Bbb H}^{3}$. If we let $q$ vary 
with velocity $\dot q$ while $y$ stays fixed, the 
corresponding derivative $\dot \Theta \left ( 
q,y\right )$ of $\Theta \left ( q,y\right )$ 
depends linearly on $\dot q$ and continuously on 
$\left ( q,y\right )$. It follows that $\dot 
\Theta \left ( q,y\right )$ is an $O\left ( \left 
\Vert \dot q\right \Vert \right )$.  
\nobreak\par\nobreak 	If we let $y$ vary with 
velocity $\dot y$ while $q$ stays fixed, the 
corresponding derivative $\dot \Theta \left ( 
q,y\right )$ depends linearly on $\dot y$ and 
differentiably on $\left ( q,y\right )$. In 
addition, if $q$ is equal to $p$, $\Theta \left ( 
q,y\right )$ is constantly equal to $\pi $ under 
such a variation, so that $\dot \Theta \left ( 
p,y\right )$ is equal to 0. The same compactness 
arguments as above now shows that $\dot \Theta 
\left ( q,y\right )$  is an $O\left ( d\left ( 
p,q\right )\left \Vert \dot y\right \Vert \right 
)$. 
\nobreak\par\nobreak 	The case of a general 
variation follows from these two cases by 
linearity of the differential of $\Theta $.  
\qed\enddemo


\nobreak\par\nobreak 	We can apply Lemma~6 to the 
tetrahedron $T$ with vertices $p=\varphi 
_{t}\left ( p_{t}\right )$, $q=\varphi _{t}\left 
( q_{t}\right )$, $y=\varphi _{t}\bigl ( 
y_{i}^{(t)}\bigr ) $ and $x$, where $x$ is the 
end point at infinity that is common to the 
geodesic containing $\varphi _{t}\left ( 
p_{t}\right )$ and $\varphi _{t}\left ( 
r_{t}\right )$ and the geodesic containing 
$\varphi _{t}\left ( q_{t}\right )$ and $\varphi 
_{t}\left ( s_{t}\right )$. In the pyramid 
$P_{t}$, $\theta _{P_{t}}\bigl ( 
p_{t}y_{i}^{(t)}\bigr )  +  \theta _{P_{t}}\bigl 
( q_{t}y_{i}^{(t)}\bigr ) $ is equal to $\Theta 
\left ( q,y\right )$ or to $2\pi -\Theta \left ( 
q,y\right )$, according to where $x$ sits with 
respect to $\varphi _{t}\left ( p_{t}\right )$ 
and $\varphi _{t}\left ( r_{t}\right )$. When we 
differentiate with respect to $t$, the variation 
$\dot p$ of $p$ is equal to 0, since $p$ is 
constant. By Lemma~5, it follows that the 
variation $\dot q$ of $q$ is an $O\left ( d\left 
( p_{t},q_{t}\right )r\left ( R_{0}\right )\right 
)$, and that the variation $\dot y$ of $y$ is an 
$O\left ( r\left ( R_{0}\right )\right )$. Since 
there is a constant $A>0$ such that $d\left ( 
p_{t},q_{t}\right )=O\bigl ( e^{-Ar\left ( 
R_{0}\right )}\bigr ) $ and since the distance 
from $y_{i}^{(t)}$ to $f_{t}\left ( \lambda 
_{t}\right )$ is bounded away from 0, Lemma~6 
shows that 
$${d\over dt} \theta _{P_{t}}\bigl ( 
p_{t}y_{i}^{(t)}\bigr )  + {d\over dt} \theta 
_{P_{t}}\bigl ( q_{t}y_{i}^{(t)}\bigr ) =O\bigl ( 
\left ( r\left ( R_{0}\right )\right )e^{-Ar\left 
( R_{0}\right )}\bigr ) .  \tag  11$$
\nobreak\par\nobreak 	The lengths $l\bigl ( 
p_{t}y_{i}^{(t)}\bigr ) $ and $l\bigl ( 
q_{t}y_{i}^{(t)}\bigr ) $ are uniformly bounded, 
and
$$l\bigl ( p_{t}y_{i}^{(t)}\bigr )  - l\bigl ( 
q_{t}y_{i}^{(t)}\bigr )  = O\left ( l\left ( 
p_{t}q_{t}\right )\right ) = O\bigl ( e^{-Ar\left 
( R_{0}\right )}\bigr ) .$$
Since $d\theta _{P_{t}}\bigl ( 
q_{t}y_{i}^{(t)}\bigr ) /dt = O\left ( r\left ( 
R_{0}\right )\right )$ by Lemma~5, it follows 
that
$$l\bigl ( p_{t}y_{i}^{(t)}\bigr )  {d\over dt} 
\theta _{P_{t}}\bigl ( p_{t}y_{i}^{(t)}\bigr )   
+ l\bigl ( q_{t}y_{i}^{(t)}\bigr )  {d\over dt} 
\theta _{P_{t}}\bigl ( q_{t}y_{i}^{(t)}\bigr ) 
 =O\bigl ( r\left ( R_{0}\right )e^{-Ar\left ( 
R_{0}\right )}\bigr ) .  \tag  12$$
Since, for every $r\geqslant 1$, the number of 
components $R_{0}$ of $R_{i}^{(0)}-\lambda _{0}$ 
such that $r\left ( R_{0}\right )=r$ is uniformly 
bounded, \rom{(12)} guarantees the convergence of 
the series
$${\textstyle {1\over 2}} \sum _{R_{t}} \left ( 
l\bigl ( p_{t}y_{i}^{(t)}\bigr )  {d\over dt} 
\theta _{P_{t}}\bigl ( p_{t}y_{i}^{(t)}\bigr )  + 
l\bigl ( q_{t}y_{i}^{(t)}\bigr )  {d\over dt} 
\theta _{P_{t}}\bigl ( q_{t}y_{i}^{(t)}\bigr )  
\right ), \tag  13$$
as we sum over all components $R_{t}$ of 
$R_{i}^{(t)}-\lambda _{t}$.  
\nobreak\par\nobreak 	The convergence of the sums
$${\textstyle {1\over 2}}  \sum _{R_{t}} \left ( 
l\bigl ( r_{t}y_{i}^{(t)}\bigr )  {d\over dt} 
\theta _{P_{t}}\bigl ( r_{t}y_{i}^{(t)}\bigr )  + 
l\bigl ( s_{t}y_{i}^{(t)}\bigr )  {d\over dt} 
\theta _{P_{t}}\bigl ( s_{t}y_{i}^{(t)}\bigr )  
\right ), \tag  14 $$
and
$${\textstyle {1\over 2}}  \sum _{R_{t}} \left ( 
l\left ( p_{t}r_{t}\right ) {d\over dt} \theta 
_{P_{t}}\bigl ( p_{t}r_{t}\bigr )  + l\left ( 
s_{t}q_{t}\right ) {d\over dt} \theta 
_{P_{t}}\bigl ( s_{t}q_{t}\bigr )  \right )  \tag 
 15$$
is proved by arguments which are, identical to 
the above one for \rom{(14)}, and very similar 
for \rom{(15)}.
\nobreak\par\nobreak 	Now, consider the 
tetrahedron $T_{t}$, with vertices $p_{t}$, 
$q_{t}$, $x_{i}^{(t)}$, $y_{i}^{(t)}$, which is 
also associated to $R_{t}$. Orient the boundary 
$\partial T_{t}$ so that the orientation induced 
on the triangle $p_{t}q_{t}y_{i}^{(t)}$ is 
opposite to the orientation induced by the 
boundary $\partial P_{t}$, and orient $T_{t}$ 
accordingly. The Schl{\"a}fli formula expresses 
$d \roman {vol}\left ( T_{t}\right ) /dt$ as a 
sum of terms corresponding to its edges. Summing 
over all components $R_{t}$ of 
$R_{i}^{(t)}-\lambda _{t}$, we obtain the 
following four sums, whose convergence is proved 
by arguments similar to the ones used for 
$P_{t}$.
$${\textstyle {1\over 2}} \sum _{R_{t}} l\left ( 
p_{t}q_{t}\right ){d\over dt}\theta _{T_{t}}\bigl 
( p_{t}q_{t}\bigr ) , \tag  16 $$
$${\textstyle {1\over 2}}  \sum _{R_{t}} l\bigl ( 
x_{i}^{(t)}y_{i}^{(t)}\bigr ) {d\over dt}\theta 
_{T_{t}}\bigl ( x_{i}^{(t)}y_{i}^{(t)}\bigr ) , 
\tag  17 $$
$${\textstyle {1\over 2}}  \sum _{R_{t}} \left ( 
l\bigl ( p_{t}x_{i}^{(t)}\bigr ) {d\over 
dt}\theta _{T_{t}}\bigl ( p_{t}x_{i}^{(t)}\bigr ) 
			+ l\bigl ( q_{t}x_{i}^{(t)}\bigr ) {d\over 
dt}\theta _{T_{t}}\bigl ( q_{t}x_{i}^{(t)}\bigr ) 
 \right ) \tag  18$$
and
$${\textstyle {1\over 2}}  \sum _{R_{t}} \left (  
l\bigl ( p_{t}y_{i}^{(t)}\bigr ) {d\over 
dt}\theta _{T_{t}}\bigl ( p_{t}y_{i}^{(t)}\bigr ) 
			+ l\bigl ( q_{t}y_{i}^{(t)}\bigr ) {d\over 
dt}\theta _{T_{t}}\bigl ( q_{t}y_{i}^{(t)}\bigr ) 
 \right ). \tag  19$$
converge, by arguments similar to the ones used 
for $P_{t}$.
\nobreak\par\nobreak 	Since the convergence of 
all these sums is uniform in $t$, we conclude 
that the volume of $H_{t}\bigl ( 
R_{i}^{(t)}\times \left \lbrack 0,1\right \rbrack 
\bigr ) $ is differentiable in $t$, and that
$$
{d\over dt}\roman {vol}\bigl ( H_{t}\bigl ( 
R_{i}^{(t)}\times \left \lbrack 0,1\right \rbrack 
\bigr ) \bigr ) ={d\over dt}\sum _{R_{t}} \roman 
{vol}\left ( P_{t}\right ) + {d\over dt}\sum 
_{R_{t}} \roman {vol}\left ( T_{t}\right ) $$
is equal to the sum of all terms \rom{(9--10)}, 
\rom{(13--15)} and \rom{(16--19)} as $R_{t}$ 
ranges over all components of 
$R_{i}^{(t)}-\lambda _{t}$. (We are here using an 
abuse of notation, where $H_{t}\bigl ( 
R_{i}^{(t)}\times \left \lbrack 0,1\right \rbrack 
\bigr ) $ represents the chain defined by 
restriction of $H_{t}$ to $R_{i}^{(t)}\times 
\left \lbrack 0,1\right \rbrack $ and not the 
image of this map. In particular, the volume is 
computed by taking into account the sign of the 
Jacobian of $H_{t}$, and may very well be 
negative. We will use the same abuse of notation 
below when considering the boundary of this 
chain.)
\nobreak\par\nobreak 	The term \rom{(17)} is 
particularly simple. Indeed, consider the corners 
$p_{i}^{(t)}=f_{t}\left ( a,c\right )$, 
$q_{i}^{(t)}=f_{t}\left ( b,c\right )$, 
$r_{i}^{(t)}=f_{t}\left ( a,d\right )$, 
$s_{i}^{(t)}=f_{t}\left ( b,d\right )$ of the 
image of the rectangle $R_{i}^{(t)}\cong \left 
\lbrack a,b\right \rbrack \times \left \lbrack 
c,d\right \rbrack $ under $f_{t}$, as in 
Figure~1. Then, the sum \rom{(17)} is equal to

$${\textstyle {1\over 2}} l\bigl ( 
x_{i}^{(t)}y_{i}^{(t)}\bigr )  {d\over dt}\theta 
_{H_{t}\bigl ( R_{i}^{(t)}\times \left \lbrack 
0,1\right \rbrack \bigr ) }\bigl ( 
x_{i}^{(t)}y_{i}^{(t)}\bigr )  \tag  20$$
where $\theta _{H_{t}\bigl ( R_{i}^{(t)}\times 
\left \lbrack 0,1\right \rbrack \bigr ) }\bigl ( 
x_{i}^{(t)}y_{i}^{(t)}\bigr )  $ is the external 
dihedral angle between the triangles 
$x_{i}^{(t)}y_{i}^{(t)}p_{i}^{(t)}$ and 
$x_{i}^{(t)}y_{i}^{(t)}q_{i}^{(t)}$. Indeed, this 
follows from the fact that
$$\pi -\theta _{H_{t}\bigl ( R_{i}^{(t)}\times 
\left \lbrack 0,1\right \rbrack \bigr ) }\bigl ( 
x_{i}^{(t)}y_{i}^{(t)}\bigr )   
			=  \sum _{R_{t}} \left ( \pi -\theta 
_{T_{t}}\bigl ( x_{i}^{(t)}y_{i}^{(t)}p_{t}, 
y_{i}^{(t)}x_{i}^{(t)}q_{t}\bigr )  \right ) \tag 
21$$
in ${\Bbb R}/2\pi {\Bbb Z}$ (because $f_{t}\left 
( \left ( \left \lbrack a,b\right \rbrack \times 
\left \lbrace c\right \rbrace \right )\cap 
\lambda _{t}\right )$ has 1--dimensional measure 
0) and that the convergence of the sum \rom{(17)} 
is uniform in $t$. 
\nobreak\par\nobreak 	It turns out that the terms 
\rom{(13)} and \rom{(19)} almost cancel out. 
Indeed, they both involve edges of the form 
$p_{t}y_{i}^{(t)}$ and $q_{t}y_{i}^{(t)}$. In 
general, the contributions to \rom{(13)} and 
\rom{(19)} of each individual edge  
$p_{t}y_{i}^{(t)}$ or $q_{t}y_{i}^{(t)}$ do not 
add up to 0. However, we will show that only four 
terms remain when we sum these contributions over 
all rectangles $R_{t}$. This will require the 
consideration of the bending cocycle of a certain 
pleated fan. 
\nobreak\par\nobreak 	Consider the closure 
$P_{i}^{(t)}\subset H_{t}\bigl ( 
R_{i}^{(t)}\times \left \lbrack 0,1\right \rbrack 
\bigr ) $ of the union of the pyramids $P_{t}$. 
It is partially bounded by the joint 
$F_{i}^{(t)}$ of the point $y_{i}^{(t)}$ and of 
the arc $f_{t}\left ( \left \lbrack a,b\right 
\rbrack \times \left \lbrace c\right \rbrace 
\right )$. This $F_{i}^{(t)}$ is a pleated fan 
with pleating locus the joint $\mu _{i}^{(t)}$ of 
$y_{i}^{(t)}$ and $f_{t}\left ( \left ( \left 
\lbrack a,b\right \rbrack \times \left \lbrace 
c\right \rbrace \right )\cap \lambda _{t}\right 
)$. We orient $F_{i}^{(t)}$ so that the boundary 
orientation it induces on $f_{t}\left ( \left ( 
\left \lbrack a,b\right \rbrack \times \left 
\lbrace c\right \rbrace \right )\cap \lambda 
_{t}\right )$ coincides with the one coming from 
the natural orientation of $\left \lbrack 
a,b\right \rbrack $. Using the methods of 
\cite{Bo3, {\S }7}, we can measure the bending of 
$F_{i}^{(t)}$ along $\mu _{i}^{(t)}$ by an ${\Bbb 
R}/2\pi {\Bbb Z}$--valued transverse cocycle 
$\beta _{t}$ for $\mu _{i}^{(t)}$; the crucial 
property here is that the curve $f_{t}\left ( 
\left \lbrack a,b\right \rbrack \times \left 
\lbrace c\right \rbrace \right )$ is rectifiable. 
We now interpret the quantity \rom{(13)} as, 
essentially, the length of the derivative of this 
bending cocycle.

\proclaim {Lemma 7 }   As we differentiate in 
$t$, the bending cocycle $\beta _{t}$ of the 
pleated fan $F_{i}^{(t)}$ admits a derivative 
$\dot \beta _{t}$, which is an ${\Bbb R}$--valued 
transverse cocycle for the pleating locus $\mu 
_{i}^{(t)}$. In addition, $\dot \beta _{t}$ has a 
well defined length $l\bigl ( \dot \beta 
_{t}\bigr ) $ in $M_{t}$, and the quantity 
\rom{(13)} is equal to
$$ {\textstyle {1\over 2}} l\bigl ( \dot \beta 
_{t}\bigr )  + {\textstyle {1\over 2}} l\bigl ( 
p_{i}^{(t)}y_{i}^{(t)}\bigr )  {d\over dt}\theta 
_{P_{i}^{(t)}}\bigl ( p_{i}^{(t)}y_{i}^{(t)}\bigr 
)  
		+ {\textstyle {1\over 2}} l\bigl ( 
q_{i}^{(t)}y_{i}^{(t)}\bigr )  {d\over dt}\theta 
_{P_{i}^{(t)}}\bigl ( q_{i}^{(t)}y_{i}^{(t)}\bigr 
)  . \tag  22$$ \endproclaim    
\demo {Proof}   For every $x\in \lambda _{t}\cap 
R_{i}^{(t)}$, consider that the tangent plane at 
$f_{t}\left ( x\right )\in M_{t}$ that is tangent 
to the geodesic arc $f_{t}\left ( x\right 
)y_{i}^{(t)}$ and to the image under $f_{t}$ of 
the leaf of $\lambda _{t}$ containing $x$; this 
plane is a Lipschitz function of $x$. 
Consequently, the expression of the bending 
cocycle given by Lemma~36 of \cite{Bo3} shows 
that, for every arc $k$ in $\left \lbrack 
a,b\right \rbrack \times \left \lbrace c\right 
\rbrace $ whose end points are disjoint from 
$\lambda _{t}$, 
$$
\beta _{t}\left ( k\right ) = \sum _{R_{t}\cap 
k\mathbin{\not =}\emptyset }   \left ( \theta 
_{P_{t}}\bigl ( p_{t}y_{i}^{(t)}\bigr )  -\pi 
+\theta _{P_{t}}\bigl ( q_{t}y_{i}^{(t)}\bigr )  
\right )
 - \theta _{P_{t}^{-}}\bigl ( 
p_{t}^{-}y_{i}^{(t)}\bigr )  +\pi -\theta 
_{P_{t}^{+}}\bigl ( q_{t}^{+}y_{i}^{(t)}\bigr ) 
 \tag  23$$
where the sum is taken over all components 
$R_{t}$ of $R_{i}^{(t)}-\lambda _{t}$ that meet 
$k$ and where, if $R_{t}^{+}$ and $R_{t}^{-}$ are 
the components of $R_{i}^{(t)}-\lambda _{t}$ that 
respectively contain the positive and negative 
end point of $k$, $P_{t}^{\pm }$, $p_{t}^{\pm }$, 
$q_{t}^{\pm }$ denote the pyramid and vertices 
associated to $f_{t}\left ( R_{t}^{\pm }\right )$ 
by the usual labelling conventions. From Lemma~4 
and \rom{(11)}, we conclude that $\beta _{t}\left 
( k\right )$ has a derivative $\dot \beta 
_{t}\left ( k\right )\in {\Bbb R}$ with respect 
to $t$, given by 
$$
\dot \beta _{t}\left ( k\right ) = \sum 
_{R_{t}\cap k\mathbin{\not =}\emptyset }   \left 
( {d\over dt}\theta _{P_{t}}\bigl ( 
p_{t}y_{i}^{(t)}\bigr )  +{d\over dt}\theta 
_{P_{t}}\bigl ( q_{t}y_{i}^{(t)}\bigr )  \right )
 - {d\over dt}\theta _{P_{t}^{-}}\bigl ( 
p_{t}^{-}y_{i}^{(t)}\bigr )  -{d\over dt}\theta 
_{P_{t}^{+}}\bigl ( q_{t}^{+}y_{i}^{(t)}\bigr ) 
 \tag  24$$

Since the finite additivity with respect to $k$ 
is immediate, this defines a transverse ${\Bbb 
R}$--valued cocycle $\dot \beta _{t}$ for the 
pleating locus $\mu _{i}^{(t)}$ of the pleated 
fan $F_{i}^{(t)}$.
\nobreak\par\nobreak 	In \cite{Bo1}, we showed 
how an ${\Bbb R}$--valued transverse cocycle for 
a geodesic lamination $\mu $ on $S$ defines a 
transverse H{\"o}lder distribution for $\mu $. 
However, this construction depended in a crucial 
way on some global properties of $\mu $. 
Associating a transverse H{\"o}lder distribution 
for $\mu _{i}^{(t)}$ to the transverse cocycle 
$\dot \beta _{t}$ is therefore not automatic. 
However, by \rom{(24)}, \rom{(11)} and Lemma~5, 
$\dot \beta _{t}\left ( k\right ) = O\left ( 
r\left ( R_{t}^{+}\right )\right )+O\left ( 
r\left ( R_{t}^{-}\right )\right )$ for every 
$k$, and $d\left ( p_{t}, q_{t}\right )=O\left ( 
e^{-Ar\left ( R_{t}\right )}\right )$ for every 
$R_{t}$, with the usual notation. Since, for 
every $r\geqslant 0$, the number of $R_{t}$ with 
$r\left ( R_{t}\right )=r$ is uniformly bounded, 
this is exactly what we need to use the 
techniques of \cite{Bo1} and associate to $\dot 
\beta _{t}$ a transverse H{\"o}lder distribution 
for $\mu _{i}^{(t)}$; compare \rom{(25)} below. 
In particular, we can integrate with respect to 
this distribution the length of the leaves of 
$\mu _{i}^{(t)}$, which defines the length 
$l\bigl ( \dot \beta _{t}\bigr ) $.
\nobreak\par\nobreak 	Theorem~11 of \cite{Bo1} 
provides an explicit expression for the 
transverse distribution $\dot \beta _{t}$, which 
gives  
$$l\bigl ( \dot \beta _{t}\bigr )  = \sum 
_{R_{t}} \dot \beta _{t}\left ( k\left ( 
R_{t}\right )\right ) \left ( l\bigl ( 
p_{t}y_{i}^{(t)}\bigr )  - l\bigl ( 
q_{t}y_{i}^{(t)}\bigr ) \right ) + \dot \beta 
_{t}\bigl ( f_{t}\left ( \left \lbrack a,b\right 
\rbrack \times \left \lbrace c\right \rbrace 
\right )\bigr )  l\bigl ( 
q_{i}^{(t)}y_{i}^{(t)}\bigr )  \tag  25$$
where $k\left ( R_{t}\right )$ denotes the arc in 
$\left \lbrack a,b\right \rbrack \times \left 
\lbrace c\right \rbrace $ that joins $\left ( 
a,c\right )$ to an arbitrary  point in the 
interior of $\left ( \left \lbrack a,b\right 
\rbrack \times \left \lbrace c\right \rbrace 
\right )\cap R_{t}$. Also, the Gap Lemma of 
\cite{Bo1} shows that 
$$   l\bigl ( q_{t}y_{i}^{(t)}\bigr )  = \sum 
_{R_{t}'\cap k\left ( R_{t}\right )=\emptyset }   
\left ( l\bigl ( p_{t}'y_{i}^{(t)}\bigr ) -l\bigl 
( q_{t}'y_{i}^{(t)}\bigr ) \right ) +  l\bigl ( 
q_{i}^{(t)}y_{i}^{(t)}\bigr )    \tag  26$$
where the sum is over those components $R_{t}'$ 
of $R_{i}^{(t)}-\lambda _{t}$ which do not meet 
the arc $k\left ( R_{t}\right )$, and where 
$p_{t}'$, $q_{t}'$, $r_{t}'$, $s_{t}'$ denote the 
vertices of $f_{t}\left ( R_{t}'\right )$ with 
the usual labelling conventions. Combining 
\rom{(24)}, \rom{(25)}, \rom{(26)} and 
rearranging terms, we conclude that

$$\multline
l\bigl ( \dot \beta _{t}\bigr )   = \sum _{R_{t}} 
\left ( l\bigl ( p_{t}y_{i}^{(t)}\bigr )  {d\over 
dt} \theta _{P_{t}}\bigl ( p_{t}y_{i}^{(t)}\bigr 
)  + l\bigl ( q_{t}y_{i}^{(t)}\bigr )  {d\over 
dt} \theta _{P_{t}}\bigl ( q_{t}y_{i}^{(t)}\bigr 
)  \right )  \\
- l\bigl ( p_{i}^{(t)}y_{i}^{(t)}\bigr )  {d\over 
dt}\theta _{P_{i}^{(t)}}\bigl ( 
p_{i}^{(t)}y_{i}^{(t)}\bigr )  
		- l\bigl ( q_{i}^{(t)}y_{i}^{(t)}\bigr )  
{d\over dt}\theta _{P_{i}^{(t)}}\bigl ( 
q_{i}^{(t)}y_{i}^{(t)}\bigr )  . 
\endmultline  \tag  27$$ 
This concludes the proof of Lemma~7. \qed\enddemo

\nobreak\par\nobreak 	Similarly, if $T_{i}^{(t)}$ 
denotes the closure of the union of the 
tetrahedra $T_{t}$, the quantity \rom{(18)} is 
equal to
$$ -{\textstyle {1\over 2}} l\bigl ( \dot \beta 
_{t}\bigr )  + {\textstyle {1\over 2}} l\bigl ( 
p_{i}^{(t)}y_{i}^{(t)}\bigr )  {d\over dt}\theta 
_{T_{i}^{(t)}}\bigl ( p_{i}^{(t)}y_{i}^{(t)}\bigr 
)  
		+ {\textstyle {1\over 2}} l\bigl ( 
q_{i}^{(t)}y_{i}^{(t)}\bigr )  {d\over dt}\theta 
_{T_{i}^{(t)}}\bigl ( q_{i}^{(t)}y_{i}^{(t)}\bigr 
)   \tag  28$$
where the negative sign comes from the fact that 
$F_{i}^{(t)}$ now occurs with the opposite 
orientation. Combining \rom{(22)} and \rom{(28)}, 
we conclude that the infinite sums \rom{(13)} and 
\rom{(19)} add up to the finite sum
$$ l\bigl ( p_{i}^{(t)}y_{i}^{(t)}\bigr )  
{d\over dt}\theta _{H_{t}\bigl ( 
R_{i}^{(t)}\times \left \lbrack 0,1\right \rbrack 
\bigr ) }\bigl ( p_{i}^{(t)}y_{i}^{(t)}\bigr )  
		+  l\bigl ( q_{i}^{(t)}y_{i}^{(t)}\bigr )  
{d\over dt}\theta _{H_{t}\bigl ( 
R_{i}^{(t)}\times \left \lbrack 0,1\right \rbrack 
\bigr ) }\bigl ( q_{i}^{(t)}y_{i}^{(t)}\bigr )   
\tag  29$$

\nobreak\par\nobreak 	Similarly, many terms 
cancel out as we take the sum of all terms 
\rom{(9)}, \rom{(10)}, \rom{(14)}, \rom{(15)}, 
\rom{(16)} and \rom{(18)} over the finitely many 
rectangles $R_{i}^{(t)}$. 
\nobreak\par\nobreak 	For the terms \rom{(9)}, 
\rom{(10)} and \rom{(16)}, this occurs in term by 
term cancellations. Indeed, with finitely many 
exceptions, the edge $p_{t}q_{t}$ of the 
rectangle $f_{t}\left ( R_{t}\right )$ 
corresponding to the component $R_{t}$ of 
$R_{i}^{(t)}-\lambda _{t}$ coincides with an edge 
$p_{t}'q_{t}'$ or $r_{t}'s_{t}'$ of a rectangle 
$f_{t}\left ( R_{t}'\right )$ corresponding to a 
component $R_{t}'$ of $R_{j}^{(t)}-\lambda _{t}$, 
for some rectangle $R_{j}^{(t)}$ possibly (and 
usually) different from $R_{i}^{(t)}$. The 
exceptions occur for those $p_{t}q_{t}$ which 
fall at the junction between two $f_{t}\bigl ( 
R_{j}^{(t)}\bigr ) $ and $f_{t}\bigl ( 
R_{k}^{(t)}\bigr ) $. If $p_{t}q_{t}$ coincides 
with such a $r_{t}'s_{t}'$, then
$$\theta _{P_{t}}\left ( p_{t}q_{t}\right 
)+\theta _{T_{t}}\left ( p_{t}q_{t}\right ) + 
\theta _{P_{t}'}\left ( r_{t}'s_{t}'\right ) = 
2\pi  .$$
It follows that the derivatives of these angles 
add up to 0, and that the contributions of 
$p_{t}q_{t}=r_{t}'s_{t}'$ to \rom{(9--10)} and 
\rom{(16)} cancel out. A similar argument holds 
when $p_{t}q_{t}$ is equal to $p_{t}'q_{t}'$, and 
when $r_{t}s_{t}$ is equal to some $p_{t}'q_{t}'$ 
or $r_{t}'s_{t}'$. It follows that, as we take 
the sum of all terms \rom{(9)}, \rom{(10)} and 
\rom{(16)} over all rectangles $R_{i}^{(t)}$, we 
are left with only finitely many boundary terms. 
\nobreak\par\nobreak 	By an argument analogous to 
Lemma~7, the sum \rom{(14)} can be interpreted in 
terms of the bending cocycle of the pleated fan 
that is the joint of $y_{i}^{(t)}$ and of the arc 
$f_{t}\left ( \left \lbrack a,b\right \rbrack 
\times \left \lbrace c\right \rbrace \right )$. A 
similar interpretation holds for the sum 
\rom{(18)}. It follows that, as we sum over all 
rectangles $R_{i}^{(t)}$, the terms \rom{(14)} 
and \rom{(18)} add up only to the sum of finitely 
many boundary terms.
\nobreak\par\nobreak 	Finally, it remains to 
consider the sum \rom{(15)}. The same arguments 
as in  Lemma~7 express \rom{(15)} as  
${\textstyle {1\over 2}}$ times the length of the 
derivative of the bending cocycle of the pleated 
rectangle $f_{t}\bigl ( R_{i}^{(t)}\bigr ) $, 
plus two boundary terms. It follows that, as we 
sum over all rectangles $R_{i}^{(t)}$, the terms 
\rom{(15)} add up to the sum of finitely many 
boundary terms and of ${\textstyle {1\over 
2}}l\bigl ( \dot \beta _{t}\bigr ) $.
\nobreak\par\nobreak 	Combining these analyses, 
we conclude that $d\roman {vol}\left ( H_{t}\left 
( R^{(t)}\times \left \lbrack 0,1\right \rbrack 
\right )\right )/dt$ is the sum of ${\textstyle 
{1\over 2}}l\bigl ( \dot \beta _{t}\bigr ) $ and 
of finitely many boundary terms corresponding to 
the edges of the polyhedral surface $H_{t}\left ( 
\partial R^{(t)}\times \left \lbrack 0,1\right 
\rbrack \right )$. 
\nobreak\par\nobreak 	The volume bounded by 
$f_{t}$ is the sum of the volume of $H_{t}\left ( 
R^{(t)}\times \left \lbrack 0,1\right \rbrack 
\right )$, the volume of $H_{t}\left ( \left ( 
S-R^{(t)}\right )\times \left \lbrack 0,1\right 
\rbrack \right )$, and the volume bounded by 
$g_{t}$. The last two of these volumes are 
bounded by polyhedral surfaces, and their 
variation is therefore given by Corollary~2, as 
the sum of finitely many boundary terms. We saw 
that the volume of $H_{t}\left ( R^{(t)}\times 
\left \lbrack 0,1\right \rbrack \right )$ is the 
sum of ${\textstyle {1\over 2}}l\bigl ( \dot 
\beta _{t}\bigr ) $ and of finitely many boundary 
terms. The boundary terms cancel out as in the 
proof of Corollary~2, and we conclude that the 
derivative of the volume enclosed by $f_{t}$ is 
equal to ${\textstyle {1\over 2}}l\bigl ( \dot 
\beta _{t}\bigr ) $.
\nobreak\par\nobreak 	This concludes the proof of 
Theorem~3 when there are no cusps.  \qed\enddemo 

\demo {Proof of Theorem~3 in the presence of 
cusps}   Each cusp of $M_{0}$ has a neighborhood 
of the form $B/\Gamma _{1}$, where $B$ is a 
horoball of ${\Bbb H}^{3}$ and where $\Gamma 
_{1}$ is a parabolic subgroup of $\Gamma $ that 
is isomorphic to ${\Bbb Z}$ (for a rank 1 cusp) 
or to ${\Bbb Z}^{2}$ (for a rank 2 cusp). In 
addition, because the pleating locus of $f_{0}$ 
is compact, we can choose these cusp 
neighborhoods so that the intersection of $f_{0}$ 
with the cusp neighborhoods consists of finitely 
many totally geodesic annuli leading to the 
cusps. The fact that $f_{0}$ bounds a finite 
volume $3$--chain implies that, in each cusp 
neighborhood $B/\Gamma _{1}$, these annuli bound 
a locally finite 3--chain relative to the 
boundary. If, in $f_{0}$, we chop off these 
annuli along piecewise geodesic simple closed 
curves, we suitably reconnect the pieces by 
polyhedral annuli, and we add a few polyhedral 
tori separating the rank 2 cusps from the rest of 
$M_{0}$, we obtain a compact surface $g_{0}$ 
which is homologous to 0 in $M_{0}$, and whose 
pleating locus consists of $\lambda $ and of 
finitely many edges. The symmetric difference 
between $f_{0}$ and $g_{0}$ gives a polyhedral 
surface $h_{0}$ which bounds a finite volume 
locally finite 3--chain, and such that the 
2--chain $f_{0}-g_{0}-h_{0}$ bounds a finite 
3--chain of volume 0. 
\nobreak\par\nobreak 	From $f_{t}$ in $M_{t}$, we 
can similarly define $g_{t}$ and $h_{t}$ in such 
a way that they depend differentiably on $t$. 
Then, the proof of Theorem~3 in the case without 
cusps immediately extends to show that the 
derivative at $t=0$ of the volume enclosed by 
$g_{t}$ is equal to 
$${\textstyle {1\over 2}}l_{0}\left ( \dot 
b_{0}\right ) + {\textstyle {1\over 2}} \sum _{e} 
l_{0}\left ( e\right )\dot \theta _{0}\left ( 
e\right ) \tag  30$$
where the sum is over the edges $e$ of the 
polyhedral part of $g_{0}$, and where $\theta 
_{t}\left ( e\right )$ is the external dihedral 
angle of $g_{t}$ at $e$. Note that every edge $e$ 
of $g_{t}$ occurs as an edge of $h_{t}$ with 
external dihedral angle $\pi -\theta _{t}\left ( 
e\right )$. We can then invoke an easy extension 
of the Schl{\"a}fli formula to polyhedral 
surfaces that have a compact set of edges and 
bound a finite volume locally finite 3--chain 
(Hint for a proof: cut off this locally finite 
extension by polyhedral surfaces that are 
arbitrarily close to the cusps, apply Theorem~2, 
and pass to the limit), which says that the 
derivative at $t=0$ of the volume enclosed by 
$h_{t}$ is equal to 
$$- {\textstyle {1\over 2}} \sum _{e} l_{0}\left 
( e\right )\dot \theta _{0}\left ( e\right ). 
\tag  31$$
Since the 2--chain $f_{t}-g_{t}-h_{t}$ bounds a 
volume 0 chain, adding up \rom{(30)} and 
\rom{(31)} completes the proof.      \qed\enddemo 

\remark {Remark \rom{1}}   A more attractive 
approach to the proof of Theorem~3 would be to 
approximate the pleated surface $f_{t}$ by 
polyhedral surfaces $f_{t}'$ and to show that, as 
the approximation gets better, the derivative 
given by the Schl{\"a}fli formula for the volume 
enclosed by $f_{t}'$ gets arbitrarily close to 
${\textstyle {1\over 2}}l_{t}\bigl ( \dot \beta 
_{t}\bigr ) $. This would decrease the cumbersome 
administration of building blocks in the above 
proof, and eliminate the consideration of 
internal edges whose contributions are eventually 
shown to cancel out. However, the author was 
unable to develop an approximation scheme where 
he could rigorously prove that this really 
happens.  \endremark

\remark {Remark \rom{2}}   Theorem~3 easily 
generalizes to the case where the pleating locus 
$\lambda $ of $f_{t}$ is non-compact. Indeed, 
there is a neighborhoood of the cusps of $S$ 
which meets only finitely many leaves of $\lambda 
$; see for instance \cite{CEG, Theorem~4.2.8}. 
The bending cocycle $b_{t}\in {\Cal H}\left ( 
\lambda ;{\Bbb R}/2\pi {\Bbb Z}\right )$ has the 
additional property that, for every cusp, the 
$b_{t}$--masses of the finitely many leaf ends of 
$\lambda $ converging to that cusp add up to 0. 
Since the same property holds for the derivative 
$\dot b_{0}\in {\Cal H}\left ( \lambda ;{\Bbb 
R}\right )$, this enables one to define a {\it 
finite\/} length $l_{0}\bigl ( \dot b_{0}\bigr ) 
$ as the contributions of the leaf ends 
converging to the cusps cancel out in the limit. 
The proof of Theorem~3 immediately generalizes to 
show that the equality $\dot V_{0}= {\textstyle 
{1\over 2}}l_{0}\bigl ( \dot b_{0}\bigr ) $ also 
holds in this case. \endremark

\head {\S }3. Proof of the main theorem\endhead   

\nobreak\par\nobreak 	We now prove the main 
theorem.

\proclaim {Theorem~8 }   Let $M_{t}$, $t\in \left 
[0,\varepsilon \right [$, be a cusp-preserving 
deformation of the geometrically finite 
hyperbolic $3$--manifold $M_{0}$. Let $b_{t}\in 
{\Cal M}{\Cal L}\left ( \partial C_{M_{0}}\right 
)$ be the bending measured geodesic lamination of 
the boundary $\partial C_{M_{t}}$ of the convex 
core of $M_{t}$ (using the convention that 
$\partial C_{M_{t}}$ is the unit normal bundle of 
$C_{M_{t}}$ when the convex core $C_{M_{t}}$ is 
$2$--dimensional). Then, the volume $V_{t}$ of 
the convex core $C_{M_{t}}$ admits a right 
derivative at $t=0$,  and
$$\dot V_{0} = {\textstyle {1\over 2}} l_{0}\bigl 
( \dot b_{0}\bigr ) $$
where $l_{0}\bigl ( \dot b_{0}\bigr ) $ is the 
length of the vector $\dot b_{0}$ tangent to the 
family of bending measured laminations $b_{t}$.  
\endproclaim    

\demo {Proof of Theorem~8 when there are no 
cusps}   Let $S$ denote the surface $\partial 
C_{M_{0}}$. Then, the bending measured geodesic 
laminations $b_{t}$ belong to the space ${\Cal 
M}{\Cal L}\left ( S\right )$ of measured geodesic 
laminations on $S$. Let $f_{t}:S\rightarrow 
M_{t}$ be the pleated surface whose image is the 
boundary $\partial C_{M_{t}}$, and let $\lambda 
_{t}$ be a pleating locus for $f_{t}$ which is 
maximal among compact geodesic laminations. Note 
that $\lambda _{t}$ contains the support of the 
bending measured geodesic lamination $b_{t}$.
\nobreak\par\nobreak 	We first prove Theorem~8 
under the additional assumption that, as $t$ 
tends to $0^{+}$, the geodesic lamination 
$\lambda _{t}$ tends to a geodesic lamination 
$\lambda $ for the Hausdorff topology. In 
particular, this is always the case when the 
support of $b_{0}$ is maximal. The fact that the 
$\lambda _{t}$ are maximal imply that $\lambda $ 
is also maximal.
\nobreak\par\nobreak 	As in \cite{Bo2}, we 
identify the tangent vector $\dot b_{0}$ to a 
compact geodesic lamination endowed with a 
certain transverse cocycle. Note that the support 
of $\dot b_{0}$ is necessarily contained in the 
Hausdorff limit $\lambda $ (see for instance 
\cite{Bo2, {\S }2}), so that $\dot b_{0}$ can be 
interpreted as a transverse cocycle for $\lambda 
$. For every $t$, let $f_{t}':S\rightarrow M_{t}$ 
be the (unique) pleated surface with pleating 
locus $\lambda $. We will prove Theorem~8 by 
comparing the volume $V_{t}$ of $C_{M_{t}}$ to 
the volume $V_{t}'$ enclosed by $f_{t}'$ in 
$M_{t}$. Note that $f_{0}'=f_{0}$, but that the 
pleating locus $\lambda _{t}$ of $f_{t}$ varies 
while the pleating locus $\lambda $ of $f_{t}'$ 
is constant.
\nobreak\par\nobreak 	We will use the Stokes 
Formula to compare the volumes respectively 
enclosed by the pleated surfaces $f_{t}$ and 
$f_{t}'$. Since the theorem is otherwise trivial  
by Theorem~1, we can assume that the $M_{t}$ are 
non-compact. Then, $H^{3}\left ( M_{t};{\Bbb 
R}\right )=0$, and there exists a differential 
2--form $\omega _{t}$ such that $d\omega _{t}$ is 
the volume form of $M_{t}$. We first show that 
the $\omega _{t}$ can be chosen to depend 
differentiably on $t$. 

\proclaim {Lemma 9}   There is a family of 
differential $2$--forms such that $d\omega _{t}$ 
coincides with the volume form of $M_{t}$ on a 
neighborhood of the convex core $C_{M_{t}}$, and 
such that $\omega _{t}$ depends differentiably on 
$t$ in the following sense: If we pull back the 
form $\omega _{t}$ on $M_{t}\cong {\Bbb 
H}^{3}/\rho _{t}\left ( \Gamma \right )$ to a 
form $\widetilde  \omega _{t}$ on ${\Bbb H}^{3}$, 
then $\widetilde  \omega _{t}$ depends 
differentiably on $t$.  \endproclaim    

\demo {Proof}   We will use a celebrated result 
of J.~Moser \cite{Mos}. By \cite{Mar, {\S }9}, 
there is for every $t$ a diffeomorphism $\varphi 
_{t}:M_{0}\rightarrow M_{t}$. In addition, the 
proof of this result makes it clear that $\varphi 
_{t}$ can be chosen to depend differentiably on 
$t$, in the sense that it lifts to a family of 
diffeomorphisms $\widetilde  \varphi _{t}:{\Bbb 
H}^{3}\rightarrow {\Bbb H}^{3}$ which depend 
differentiably on $t$.
\nobreak\par\nobreak 	If $\nu _{t}$ denotes the 
volume form of $M_{t}$, the $\varphi 
_{t}^{*}\left ( \nu _{t}\right )$ give a 
1--parameter family of volume forms on $M_{0}$ 
which are all cohomologous (to 0). Then, Moser's 
Lemma \cite{Mos} asserts that these volume forms 
are all isotopic: For every compact $K\subset 
M_{0}$, there are diffeomorphisms $\psi 
_{t}:M_{0}\rightarrow M_{0}$ depending 
differentiably on $t$ such that $\varphi 
_{t}^{*}\left ( \nu _{t}\right )=\psi 
_{t}^{*}\left ( \nu _{0}\right )$ on $K$ and such 
that  $\psi _{0}=\roman {Id}$ (Moser's proof 
provides a vector field, and one needs to 
restrict to a compact subset $K$ to integrate 
it). 
\nobreak\par\nobreak 	Pick a 2--form $\omega 
_{0}$ on $M_{0}$ such that $d\omega _{0}=\nu 
_{0}$, and a compact $K\subset M_{0}$ that 
contains a neighborhood of all the $\varphi 
_{t}^{-1}\left ( C_{M_{t}}\right )$. Then $\omega 
_{t}=\left ( \varphi _{t}^{-1}\right )^{*}\psi 
_{t}^{*}\left ( \omega _{0}\right )$ satisfies 
the properties required.  \qed\enddemo

\nobreak\par\nobreak 	For every $t$, the image of 
$f_{t}'$ is contained in $C_{M_{t}}$ and, because 
the complement $M-\partial C_{M_{t}}$ is 
homeomorphic to a product $S\times {\Bbb R}$ (see 
for instance \cite{EpM, {\S }1}),  $f_{t}$ is 
homologous to 0 in $C_{M_{t}}$. The Stokes 
Formula then shows that the volumes respectively 
enclosed by $f_{t}$ and $f_{t}'$ are equal to
$$V_{t}= \displaystyle \int _{S} f_{t}^{*}\left ( 
\omega _{t}\right ) \roman {\enskip and\enskip } 
V_{t}' =\displaystyle \int _{S} \left ( 
f_{t}'\right )^{*}\left ( \omega _{t}\right ). 
\tag  32$$
(We let the reader check, for instance through an 
approximation by polyhedral surfaces, that the 
Stokes formula holds for pleated surfaces). 
\nobreak\par\nobreak 	By a suitable partition of 
unity, $\omega _{t}$ coincides on a neighborhood 
of $C_{M_{t}}$ with a finite sum $\sum 
_{i=1}^{n}\omega _{t}^{(i)}$, where each $\omega 
_{t}^{(i)}$ is the push forward of a compactly 
supported 2--form $\widetilde  \omega _{t}^{(i)}$ 
on ${\Bbb H}^{3}$. In addition, using the 
diffeomorphisms $\varphi _{t}:M_{0}\rightarrow 
M_{t}$, we can arrange that the $\widetilde  
\omega _{t}^{(i)}$ depend differentiably on $t$. 
Consider the covering $\widehat S\rightarrow S$, 
pull back of the covering ${\Bbb 
H}^{3}\rightarrow M_{0}$ by the map 
$f_{0}:S\rightarrow M_{0}$, and the canonical 
lift $\widehat f_{0}:\widehat S\rightarrow {\Bbb 
H}^{3}$. Recall that $\widehat S$ consists of all 
pairs $\left ( x,y\right )\in S\times {\Bbb 
H}^{3}$ such that $f_{0}\left ( x\right )$ is 
equal to the projection of $y$ in $M_{0}$. In 
particular, the homotopy from $f_{0}$ to $f_{t}$ 
and $f'_{t}$ uniquely defines lifts $\widehat 
f_{t}$, $\widehat f'_{t}:\widehat S\rightarrow 
{\Bbb H}^{3}$. Note that, because the group $\rho 
_{t}\left ( \pi _{1}\left ( S\right )\right )$ 
acts properly discontinuously on ${\Bbb H}^{3}$, 
the pleated surfaces $\widehat f_{t}$ and 
$\widehat f_{t}'$ are proper. The definitions are 
specially designed so that

$$V_{t}=\sum _{i=1}^{n} \displaystyle \int 
_{\widehat S} \widehat f_{t}^{*}\bigl ( 
\widetilde  \omega _{t}^{(i)}\bigr )  
\roman {\enskip and\enskip } V_{t}'=\sum 
_{i=1}^{n} \displaystyle \int _{\widehat S} \bigl 
( \widehat f_{t}'\bigr ) ^{*}\bigl ( \widetilde  
\omega _{t}^{(i)}\bigr )  .  \tag  33$$
\proclaim {Lemma~10 }   For every compactly 
supported differential $2$--form $\widetilde  
\omega $ on ${\Bbb H}^{3}$, the following two 
right derivatives exist and are equal:
$${d\over dt^{+}}  \displaystyle \int _{\widehat 
S} \widehat f_{t}^{*}\bigl ( \widetilde  \omega 
\bigr )  _{ \vert t=0} 
				= {d\over dt^{+}} \displaystyle \int 
_{\widehat S} \bigl ( \widehat f_{t}'\bigr ) 
^{*}\bigl ( \widetilde  \omega \bigr )   _{ \vert 
t=0} . \tag 34$$
 \endproclaim    

\demo {Proof}   It clearly suffices to restrict 
attention to each component $\widehat S_{1}$ of 
$\widehat S$, projecting to a component $S_{1}$ 
of $S$. Because $f_{t}$ and $f_{t}'$ depend 
continuously on $t$, we can use a partition of 
unity to assume, without loss of generality, that 
there is a compact subset $\widetilde  B$ of the 
universal covering $\widetilde  S$ of $\widehat 
S_{1}$ (and $S_{1}$) such that the projection 
$\widetilde  S\rightarrow \widehat S_{1}$ is 
injective on $\widetilde  B$, and such that the 
intersection of the support of $\widetilde  
\omega $ with each $\widehat f_{t}\bigl ( 
\widehat S_{1}\bigr ) $ or $\widehat f_{t}'\bigl 
( \widehat S_{1}\bigr ) $ is contained in 
$\widetilde  f_{t}\bigl ( \widetilde  B\bigr ) $ 
or $\widetilde  f_{t}'\bigl ( \widetilde  B\bigr 
) $, respectively, where $\widetilde  f_{t}$ and 
$\widetilde  f_{t}'$ denote the composition of 
the projection $\widetilde  S\rightarrow \widehat 
S_{1}$ with the restrictions of $\widehat f_{t}$ 
and $\widehat f_{t}'$ to $\widehat S_{1}$. We now 
have 
$$  \displaystyle \int _{\widehat S_{1}} \widehat 
f_{t}^{*}\bigl ( \widetilde  \omega \bigr )   =   
\displaystyle \int _{\widetilde  B} \widetilde  
f_{t}^{*}\bigl ( \widetilde  \omega \bigr )   
\roman {\enskip and\enskip }
			  \displaystyle \int _{\widehat S_{1}} \bigl ( 
\widehat f_{t}'\bigr ) ^{*}\bigl ( \widetilde  
\omega \bigr )    	=  \displaystyle \int 
_{\widetilde  B} \bigl ( \widetilde  f_{t}'\bigr 
) ^{*}\bigl ( \widetilde  \omega \bigr )    . 
\tag 35$$

Then, the property of Lemma~10 is essentially 
proved in \cite{Bo4, {\S }2}, where we compare 
the two pleated surfaces $\widetilde  f_{t}$, 
$\widetilde  f_{t}':\widetilde  S\rightarrow 
{\Bbb H}^{3}$. However, minor adjustments are 
necessary because the pull back metrics $m_{t}$ 
and $m_{t}'$ induced on $S_{1}$ by $f_{t}$ and 
$f_{t}'$ may be different. 
\nobreak\par\nobreak 	Since $\widetilde  f_{t}$ 
and $\widetilde  f_{t}'$ are equivariant with 
respect to the same representation $\rho 
_{t}:\Gamma \rightarrow \roman {Isom}^{+}\left ( 
{\Bbb H}^{3}\right )$, it follows from \cite{Bo4, 
Proposition~5} that $\dot m_{0}=\dot m_{0}'$, 
while $m_{0}=m_{0}'$ since $f_{0}=f_{0}'$. Also, 
if $b_{t}'\in {\Cal H}\left ( \lambda ;{\Bbb 
R}/2\pi {\Bbb Z}\right )$ is the transverse 
cocycle describing the bending of $f_{t}'$, 
\cite{Bo4, Proposition~5} also shows that $\dot 
b_{0}=\dot b_{0}'\in {\Cal H}\left ( \lambda 
;{\Bbb R}\right )$. 
\nobreak\par\nobreak 	In \cite{Bo3} (see also 
\cite{EpM, {\S }3}), we show how to reconstruct 
the image of $\widetilde  f_{t}$ from the bending 
measured geodesic lamination $b_{t}$, considered 
as a transverse cocycle $b_{t}\in {\Cal H}\left ( 
\lambda _{t};{\Bbb R}\right )$, and from the 
shearing cocycle $s_{t}\in {\Cal H}\left ( 
\lambda _{t};{\Bbb R}\right )$ corresponding to 
$m_{t}$. More precisely, we start from the 
(un-)pleated surface $\widetilde  
g_{0}:\widetilde  S\rightarrow {\Bbb H}^{3}$ with 
pull back metric $m_{0}$ and bending measure 0, 
and we realize the geodesic laminations $\lambda 
$ and $\lambda _{t}$ by their $m_{0}$--geodesic 
representatives in $S$. Let $\widetilde  \lambda 
_{t}$ denote the preimage of $\lambda _{t}\cap 
S_{1}$ in $\widetilde  S$. Given an oriented 
geodesic $g$ of ${\Bbb H}^{3}$ and $z\in {\Bbb 
C}/2\pi i{\Bbb Z}$, let $U_{g}^{z}:{\Bbb 
H}^{3}\rightarrow {\Bbb H}^{3}$ denote the 
composition of a hyperbolic translation of signed 
length $\roman {Re}z$ along $g$ and of a 
hyperbolic rotation of angle $\roman {Im}z$ 
around $g$. Then, for every component $P$ of 
$\widetilde  S-\widetilde  \lambda _{t}$, the 
ideal triangle $\widetilde  f_{t}\left ( P\right 
)\subset {\Bbb H}^{3}$ is defined as a limit of 
terms
$$U_{\widetilde  g_{0}\left ( g_{1}^{t}\right 
)}^{c_{t}\left ( \gamma _{1}\right )} 
U_{\widetilde  g_{0}\left ( g_{2}^{t}\right 
)}^{c_{t}\left ( \gamma _{2}\right )} \dots 
U_{\widetilde  g_{0}\left ( g_{p}^{t}\right 
)}^{c_{t}\left ( \gamma _{p}\right )} \widetilde  
g_{0}\left ( P\right ) ,$$
where the $g_{i}^{t}$ are suitably chosen 
$m_{0}$--geodesics in $\widetilde  S$ and where 
the coefficients $c_{t}\left ( \gamma _{i}\right 
)\in {\Bbb C}/2\pi i{\Bbb Z}$ have their real 
part determined by the shearing cocycle $s_{t}$ 
and their imaginary part determined by the 
bending cocycle $b_{t}$. Compare \cite{Bo4, {\S 
}2}.
\nobreak\par\nobreak 	If $s_{t}'$ is the shearing 
cocycle corresponding to the metric $m_{t}'$, 
replacing $\lambda _{t}$, $s_{t}\in {\Cal H}\left 
( \lambda _{t};{\Bbb R}\right )$, $b_{t}\in {\Cal 
H}\left ( \lambda _{t};{\Bbb R}/2\pi {\Bbb 
Z}\right )$ by $\lambda $, $s_{t}'\in {\Cal 
H}\left ( \lambda ;{\Bbb R}\right )$ , $b_{t}'\in 
{\Cal H}\left ( \lambda ;{\Bbb R}/2\pi {\Bbb 
Z}\right )$ similarly defines a pleated surface 
$\widetilde  g_{t}'$ which coincides with 
$\widetilde  f_{t}$ on a neighborhood of a base 
point of $\widetilde  S$, and whose image is 
almost the same as the image of $\widetilde  
f_{t}'$. However, there is a minor subtlety here, 
in the sense that we can only conclude that there 
is an isometry $A_{t}$ of ${\Bbb H}^{3}$ such 
that $\widetilde  f_{t}'\left ( P\right 
)=A_{t}\circ \widetilde  g_{t}'\left ( P\right )$ 
for every component $P$ of $\widetilde  
S-\widetilde  \lambda $. However, $\widetilde  
g_{t}'$ is equivariant with respect to the 
representation $\rho _{t}'=A_{t}\rho 
_{t}A_{t}^{-1}:\Gamma \rightarrow \roman 
{Isom}^{+}\left ( {\Bbb H}^{3}\right )$ and 
\cite{Bo4, {\S }2} shows that $\rho '_{0}=\rho 
_{0}$ and $\dot \rho _{0}'=\dot \rho _{0}$. 
Looking at fixed points, for instance, one easily 
concludes that $A_{0}=\roman {Id}$ and $\dot 
A_{0}=0$. 
\nobreak\par\nobreak 	In \cite{Bo4, {\S }2}, we 
show that $\widetilde  f_{t}$ and $\widetilde  
g_{t}'$ are infinitesimally close as $t$ tends to 
0, and this uniformly on compact subsets of 
$\widetilde  S-\widetilde  \lambda $. In 
particular,  the arguments of \cite{Bo4, {\S }2}, 
and most notably the key growth estimate of 
\cite{Bo4, Lemma~7}, show that
$${d\over dt^{+}}  \displaystyle \int 
_{\widetilde  B} \widetilde  f_{t}^{*}\bigl ( 
\widetilde  \omega \bigr )  _{ \vert t=0} 
				= {d\over dt^{+}} \displaystyle \int 
_{\widetilde  B} \bigl ( \widetilde  g_{t}'\bigr 
) ^{*}\bigl ( \widetilde  \omega \bigr )   _{ 
\vert t=0} . \tag  36$$
because $s_{0}=s_{0}'$, $\dot s_{0}=\dot s_{0}'$, 
$b_{0}=b_{0}'$, $\dot b_{0}=\dot b_{0}'$. By 
construction,
$$ \displaystyle \int _{\widetilde  B} \bigl ( 
\widetilde  f_{t}'\bigr ) ^{*}\bigl ( \widetilde  
\omega \bigr )  =  \displaystyle \int 
_{\widetilde  B} \bigl ( \widetilde  g_{t}'\bigr 
) ^{*}A_{t}^{*}\bigl ( \widetilde  \omega \bigr ) 
  \tag  37$$
and \rom{(34)} immediately follows from 
\rom{(36-37)} and from the fact that  
$A_{0}=\roman {Id}$ and $\dot A_{0}=0$.   
\qed\enddemo 

\nobreak\par\nobreak 	From \rom{(33)} and 
\rom{(34)}, we conclude that 
$$
\dot V_{0} = \sum _{i=1}^{n} {d\over dt^{+}}  
\displaystyle \int _{\widehat S} \widehat 
f_{t}^{*}\bigl ( \widetilde  \omega 
_{0}^{(i)}\bigr )  _{ \vert t=0} 
		+  \sum _{i=1}^{n}   \displaystyle \int 
_{\widehat S} \widehat f_{0}^{*}\bigl ( {d\over 
dt^{+}}\widetilde  \omega _{t}^{(i)} {}_{ \vert 
t=0} \bigr )   \tag 38$$
where the chain rule is justified by the fact 
that the map $t\mapsto \displaystyle \int 
_{\widehat S} \widehat f_{t}^{*}\left ( 
\widetilde  \omega \right )$ is continuous, 
uniformly on compact sets for the 2--form 
$\widetilde  \omega $. The last term of 
\rom{(38)} is equal to  
$$  \displaystyle \int _{S} f_{0}^{*}\bigl (   
\sum _{i=1}^{n}  {d\over dt^{+}}\omega _{t}^{(i)} 
{}_{ \vert t=0} \bigr )  
 =   \displaystyle \int _{S} f_{0}^{*}\bigl ( 
\dot \omega _{0}\bigr )   = \displaystyle \int 
_{C_{M_{0}}} d\dot \omega _{0} = 0  \tag 39$$
since the lift of $d\omega _{t}$ to ${\Bbb 
H}^{3}$ is constant, equal to the volume form of 
${\Bbb H}^{3}$. Therefore,
$$ \dot V_{0} = \sum _{i=1}^{n} {d\over dt^{+}}  
\displaystyle \int _{\widehat S} \widehat 
f_{t}^{*}\bigl ( \widetilde  \omega 
_{0}^{(i)}\bigr )  _{ \vert t=0}   \tag  40$$
and similarly
$$ \dot V_{0} '= \sum _{i=1}^{n} {d\over dt^{+}}  
\displaystyle \int _{\widehat S} \bigl ( \widehat 
f_{t}'\bigr ) ^{*}\bigl ( \widetilde  \omega 
_{0}^{(i)}\bigr )  _{ \vert t=0}  . \tag  41$$
\nobreak\par\nobreak 	Combining \rom{(40)}, 
\rom{(41)}, Lemma~10 and Theorem~3, 
$$ \dot V_{0} = \dot V_{0} '= {\textstyle {1\over 
2}}l_{0}\bigl ( \dot b_{0}'\bigr )  =  
{\textstyle {1\over 2}}l_{0}\bigl ( \dot 
b_{0}\bigr )  $$
which concludes the proof of Theorem~8 under the 
assumption that, as $t$ tends to $0^{+}$, the 
geodesic lamination $\lambda _{t}$ converge to a 
geodesic lamination $\lambda $ for the Hausdorff 
topology. 
\nobreak\par\nobreak 	In the general case, choose 
a sequence $t_{n}$ converging to 0 such that 
$\lambda _{t_{n}}$ converges to some geodesic 
lamination $\lambda $. Then, the arguments of the 
particular case apply to show that $\left ( 
V_{t_{n}}-V_{0}\right )/t_{n}$ tends to $  
{\textstyle {1\over 2}}l_{0}\bigl ( \dot 
b_{0}\bigr )  $ as $n$ tends to infinity. Since 
this holds for any such sequence $t_{n}$, we 
conclude that $  \dot V_{0} = {\textstyle {1\over 
2}}l_{0}\bigl ( \dot b_{0}\bigr )  $ in the 
general case as well.   \qed\enddemo

\demo {Proof of Theorem~8 in the presence of 
cusps}    When the manifolds $M_{t}$ have cusps, 
the boundary $\partial C_{M_{t}}$ is totally 
geodesic near the cusps of $M_{t}$. Therefore, we 
can use the same technique as in the proof of 
Theorem~3, and chop  off pieces of $C_{M_{t}}$ by 
polyhedral surfaces near the cusps. The proof in 
this case then follows from the proof in the case 
without cusps, as for Theorem~3. \qed\enddemo 

\head {\S }4. Convex cores with totally geodesic 
boundary\endhead    

\nobreak\par\nobreak 	We conclude with an easy 
application of Theorem~8. 

\proclaim {Corollary 11 }   Let $M$ be a 
geometrically finite hyperbolic $3$--manifold 
whose convex core $C_{M}$ has totally geodesic 
boundary, but is not $2$--dimensional. Consider 
the volumes of the convex cores of the 
cusp-preserving deformations of $M$. Then $M$ is 
a strict local minimum for this volume function. 
\endproclaim    
\demo {Proof}   Let ${\Cal Q}{\Cal D}\left ( 
M\right )$ denote the space of hyperbolic 
3--manifolds obtained by cusp-preserving 
deformations of $M$. Theorem~8 determines the 
tangent map $T_{M}V:T_{M}{\Cal Q}{\Cal D}\left ( 
M\right )\rightarrow {\Bbb R}$ of the function 
$V:{\Cal Q}{\Cal D}\left ( M\right )\rightarrow 
{\Bbb R}^{+}$ defined by consideration of the 
volumes of convex cores, in terms of the tangent 
map of the bending measured lamination map $\beta 
:{\Cal Q}{\Cal D}\left ( M\right )\rightarrow 
{\Cal M}{\Cal L}\left ( \partial C_{M}\right )$ 
analyzed in \cite{Bo4}. Note that these tangent 
maps are not necessarily linear; see \cite{Bo4, 
{\S }1}.
\nobreak\par\nobreak 	Suppose that $M$ is not a 
strict local minimum for the volume function $V$. 
Then, from a sequence $M_{n}\in {\Cal Q}{\Cal 
D}\left ( M\right )$ converging to $M$ with 
$V\left ( M_{n}\right )\leqslant V\left ( M\right 
)$, we can construct a non-zero tangent vector 
$v\in T_{M}{\Cal Q}{\Cal D}\left ( M\right )$ 
such that $T_{M}V\left ( v\right )\leqslant 0$. 
Now, Theorem~11 says that $T_{M}V\left ( v\right 
) = {\textstyle {1\over 2}} l\left ( T_{M}\beta 
\left ( v\right )\right )$. Since the boundary of 
$C_{M}$ is totally geodesic, $\beta \left ( 
M\right )=0$ and $T_{M}\beta \left ( v\right )$ 
is a vector tangent to ${\Cal M}{\Cal L}\left ( 
S\right )$ at 0. By \cite{Bo2, Theorem~21}, 
$T_{M}\beta \left ( v\right )$ is therefore a 
geodesic lamination with a transverse (positive) 
measure. The main consequence of this is that 
$T_{M}\beta \left ( v\right )$ has positive 
length if $T_{M}\beta \left ( v\right 
)\mathbin{\not =}0$.  Therefore, since 
$T_{M}V\left ( v\right )\leqslant 0$, we must 
have $T_{M}\beta \left ( v\right )=0$. The proof 
is then completed by the following lemma.

\proclaim {Lemma~12 }   Under the hypotheses of 
Corollary~11, there is no non-zero tangent vector 
$v\in T_{M}{\Cal Q}{\Cal D}\left ( M\right )$ 
such that $T_{M}\beta \left ( v\right )=0$. In 
other words, there is no infinitesimal 
deformation of $M$ which infinitesimally keeps 
the boundary of $C_{M}$ flat.   \endproclaim    
\demo {Proof}   Suppose there is such a tangent 
vector $v$. Consider the manifold $DM$ obtained 
by taking the double of $C_{M}$ along its 
boundary. Namely, $DM$ is obtained by gluing two 
copies of $C_{M}$ along their boundaries by the 
natural identification.  Then, because $\partial 
C_{M}$ is totally geodesic, the hyperbolic metric 
of $M$ gives a finite volume complete hyperbolic 
metric on $DM$.  
\nobreak\par\nobreak 	First suppose that, in 
addition, there is a 1--parameter family of 
deformations $M_{t}$, $t\in \left [0,\varepsilon 
\right [$, such that $M_{0}=M$, $\dot M_{0}=v$ 
and  the convex cores $C_{M_{t}}$ all have 
totally geodesic boundary. Then, the hyperbolic 
manifolds $DM_{t}$ give a non-trivial 
cusp-preserving deformation of $DM$, which is 
excluded by Mostow's Rigidity Theorem \cite{Mow}. 
\nobreak\par\nobreak 	In general, the fact that 
$T_{M}\beta \left ( v\right )=0$ only means that, 
for a 1--parameter family of deformations 
$M_{t}$, $t\in \left [0,\varepsilon \right [$, 
with $M_{0}=M$ and $\dot M_{0}=v$, the bending 
measured geodesic lamination $b_{t}\in {\Cal 
M}{\Cal L}\left ( \partial C_{M}\right )$ of 
$M_{t}$ is such that $b_{0}=0$ and $\dot 
b_{0}=0$. We will use an infinitesimal version of 
the above argument, based on the Calabi-Weil 
Rigidity Theorem \cite{Cal}\cite{Wei} which is an 
infinitesimal version of Mostow rigidity. Indeed, 
interpreting ${\Cal Q}{\Cal D}\left ( M\right )$ 
as a space of representations $\rho :\pi 
_{1}\left ( M\right )\rightarrow \roman 
{Isom}^{+}\left ( {\Bbb H}^{3}\right )$, the Weil 
machinery (see for instance \cite{Rag}) expresses 
the tangent space of ${\Cal Q}{\Cal D}\left ( 
M\right )$ at $M$ as a subspace of the cohomology 
group $H^{1}\left ( \pi _{1}\left ( M\right ), 
\roman {Ad}\right )$, where $\roman {Ad}$ denotes 
the adjoint representation of $\pi _{1}\left ( 
M\right )$ in the Lie algebra of $\roman 
{Isom}^{+}\left ( {\Bbb H}^{3}\right )$ defined 
by the holonomy of $M$. The reason why 
$T_{M}{\Cal Q}{\Cal D}\left ( M\right )$ is only 
a subspace of $H^{1}\left ( \pi _{1}\left ( 
M\right ), \roman {Ad}\right )$ is that we 
restrict attention to cusp-preserving 
deformations. If $S$ is a component of $\partial 
C_{M}$, let $M_{S}$ be the covering of $M$ with 
$\pi _{1}\left ( M_{S}\right )=\pi _{1}\left ( 
S\right )$. The metric of $M$ lifts to a Fuchsian 
hyperbolic metric on $M_{S}$. Because $T_{M}\beta 
\left ( v\right )=0$, \cite{Bo4, Proposition~5} 
shows that the differential of the restriction 
map ${\Cal Q}{\Cal D}\left ( M\right )\rightarrow 
{\Cal Q}{\Cal D}\left ( M_{S}\right )$ sends $v$ 
to a vector $v_{S}\in T_{M_{S}}{\Cal Q}{\Cal 
D}\left ( M_{S}\right )=H^{1}\left ( \pi 
_{1}\left ( S\right ), \roman {Ad}\right )$ that 
is tangent to the submanifold of Fuchsian 
deformations of $M_{S}$. In particular, $v_{S}$ 
is invariant under the automorphism of 
$T_{M_{S}}{\Cal Q}{\Cal D}\left ( M_{S}\right )$ 
induced by the isometry that reflects $M_{S}$ 
across the totally geodesic surface $C_{M_{S}}$. 
Since this holds for every component $S$ of 
$\partial C_{M}$, a Mayer-Vietoris type argument 
shows that $v$ provides a non-trivial element of 
$H^{1}\left ( \pi _{1}\left ( DM\right ), \roman 
{Ad}\right )$ which keeps the cusps parabolic. 
However, the Calabi-Weil Rigidity Theorem 
\cite{Cal}\cite{Wei}, improved by Garland 
\cite{Gar} for the case with cusps, says that 
there is no such non-trivial element of 
$H^{1}\left ( \pi _{1}\left ( DM\right ), \roman 
{Ad}\right )$. \qed\enddemo 

\nobreak\par\nobreak 	This concludes the proof of 
Corollary~11.  \qed\enddemo

\nobreak\par\nobreak 	When $C_{M}$ is 
2--dimensional, $M$ is of course a global minimum 
for $V$ since $V\left ( M\right )=0$. In this 
case, $M$ is Fuchsian or twisted Fuchsian, and 
there are many deformations which keep the convex 
core 2--dimensional. Therefore, $M$ is only a 
weak local minimum.
\nobreak\par\nobreak 	There presumably is a 
converse to Corollary~11: If $M$ is a local 
minimum for the convex core volume function $V$, 
then the boundary $\partial C_{M}$ is totally 
geodesic. This would follow from Theorem~8 and a 
conjectural extension to convex cores of Cauchy's 
Rigidity Theorem for polyhedra; compare 
\cite{RiH}.


\Refs\refstyle{A}\widestnumber\key{CEG}

\ref\key  AVS  \by D.V. Alekseevskij, E.B. 
Vinberg, A.S. Solodovnikov  \paper  Geometry of 
spaces of constant curvature\inbook Geometry II 
\ed E.B. Vinberg \bookinfo  Encyclopaedia of 
Mathematical Sciences vol. 29 \yr 1993 \publ 
Springer-Verlag \publaddr Berlin Heidelberg New 
York \pages 1--138 \endref

\ref\key BiS \by J.S.~Birman, C.~Series \paper  
Geodesics with bounded intersection numbers on 
surfaces are sparsely distributed  \jour Topology 
\vol 24  \yr 1985 \pages 217-225 \endref

\ref\key Bo1   \by F. Bonahon \paper Transverse 
H{\"o}lder distributions for geodesic laminations 
 \jour Topology \vol 36\yr 1997\pages 103--122 
\endref

\ref\key Bo2   \by F. Bonahon \paper Geodesic 
laminations with transverse H{\"o}lder 
distributions  \jour Ann. Scient. \'Ec. Norm. 
Sup. \vol 30\yr 1997 \pages 205--240 \endref

\ref\key   Bo3   \by F. Bonahon \paper Shearing 
hyperbolic surfaces, bending pleated surfaces, 
and the Thurston symplectic form  \jour Ann. Fac. 
Sci. Toulouse \vol 5 \yr 1996  \pages 
233--297\endref

\ref\key   Bo4   \by F. Bonahon \paper Variations 
of the boundary geometry of 3--dimensional 
hyperbolic convex cores  \paperinfo preprint \yr 
1996 \endref

\ref\key   Cal   \by E. Calabi \paper On compact 
riemannian manifolds with constant curvature, 
\rom{I}  \inbook Proc.
Sympos. Pure Math., vol. III   \yr 1961 \publ  
Amer. Math. Soc. \publaddr Providence\pages 
155--180\endref

\ref\key  CEG  \by R.D.~Canary, D.B.A.~Epstein, 
P.~Green  \paper  Notes on notes of 
Thurston\inbook Analytical and Geometrical 
aspects of Hyperbolic space \ed D.B.A.~Epstein 
\bookinfo  L.M.S. Lecture Notes Series vol. 111 
\yr 1987 \publ Cambridge Univ. Press \pages 3--92 
\endref

\ref\key CaB \by A. Casson, S.A. Bleiler \book 
Automorphisms of surfaces after Thurston and 
Nielsen \bookinfo  \yr 1988 \publ Cambridge 
University Press \endref

\ref\key EpM \by D.B.A. Epstein, A. Marden \paper 
Convex hulls in hyperbolic spaces, a theorem of 
Sullivan, and measured pleated surfaces \inbook 
Analytical and geometric aspects of hyperbolic 
space \ed D.B.A.~Epstein \bookinfo L.M.S. Lecture 
Note Series vol. 111 \yr 1986 \publ Cambridge 
University Press  \pages 113--253 \endref

\ref\key   Gar   \by H. Garland \paper  A 
rigidity theorem for discrete subgroups \jour 
Trans Amer. Math. Soc.  \vol 129 \yr 1967 \pages 
1--25\endref

\ref\key Kne \by H. Kneser \paper Der 
Simplexinhalt in der nichteuklidischen Geometrie 
\jour Deutsche Math. \vol 1  \yr 1936 \pages 
337--340 \endref

\ref\key Mar \by A. Marden \paper The geometry of 
finitely generated Kleinian groups \paperinfo  
\jour Ann. Math. \vol 99  \yr 1974 \pages 
383--462 \miscnote  \finalinfo  \endref

\ref\key Mos \by J. Moser \paper On the volume 
elements on a manifold\paperinfo  \jour Trans. 
Amer. Math. Soc. \vol 120 \yr 1965 \pages  
286--294 \endref

\ref\key Mow \by G.D. Mostow \book  Strong 
rigidity of locally symmetric spaces\bookinfo 
Ann. Math. Studies vol. 78 \yr 1973\publ 
Princeton Univ. Press  \endref

\ref\key PeH \by R.C. Penner, J.L. Harer \book 
Combinatorics of train tracks \bookinfo Ann. 
Math. Studies vol. 125 \yr 1992 \publ Princeton 
Univ. Press  \endref

\ref\key Rag \by M.S. Ragunathan \book  Discrete 
subgroups of Lie groups\bookinfo Ergebn. der 
Mathematik vol. 58 \yr 1972  \publ 
Springer-Verlag \publaddr Berlin Heidelberg New 
York  \endref

\ref\key RiH \by I. Rivin, C.D. Hodgson \paper A 
characterization of convex polyhedra in 
hyperbolic 3--space  \jour Invent. Math. \vol 111 
\yr 1993 \pages  77--111\moreref Corrigendum\jour 
Invent. Math. \vol 117 \yr 1994 \page 359 \endref

\ref\key Sc1 \by L. Schl{\"a}fli \paper On the 
multiple integral $\iint \dots \int  dx\thinspace 
dy\dots dz$ whose limits are 
$p_{1}=a_{1}x+b_{1}y+\dots +h_{1}z>0$, $p_{2}>0$, 
\dots , $p_{n}>0$ and $x^{2}+y^{2}+\dots 
+z^{2}=1$ \jour Q.J. Math \vol 2  \yr 1858 \pages 
269--300 \moreref{\it idem\/}\jour Q.J. Math \vol 
3  \yr 1860 \pages 54--68 \moreref {\it 
idem\/}\jour Q.J. Math \vol 3  \yr 1860 \pages 
97--108\finalinfo see also \cite{Sc2, pp. 
219--270 }  \endref

\ref\key Sc2 \by L. Schl{\"a}fli \book Gesammelte 
Mathematische Abhandlungen \bookinfo Band II   
\yr 1950\publ Birkh{\"a}user \publaddr Basel 
\endref

\ref\key  Thu   \by W.P. Thurston \book The 
topology and geometry of $3$--manifolds \yr 
1976--79 \publ Princeton University \bookinfo 
Lecture notes \endref

\ref\key   Wei   \by A. Weil \paper  On discrete 
subgroups of Lie groups \rom{II} \jour  Ann. of 
Math. \vol 75 \yr 1962 \pages 578--602\endref


\endRefs
\enddocument